\documentclass{SCAEOL}
\numberwithin{equation}{section}
\usepackage{epstopdf}
\usepackage{subfigure}
\newtheorem{alg}{\bf Algorithm}
\begin{document}

\Year{2015} %
\Month{June}
\Vol{56} %
\No{1} %
\BeginPage{1} %
\EndPage{16} %
\AuthorMark{GONG FangHui {\it et al.}}
\ReceivedDay{June 24, 2015}
\AcceptedDay{June 24, 2015}

\title{A new shift strategy for the implicitly restarted generalized second-order Arnoldi method}{}


\author{GONG FangHui}{}
\author{SUN YuQuan}{Corresponding author}

\address{School of Mathematics and Systems Science, BeiHang University, Beijing {\rm 100191}, China;}
\Emails{ fanghui.gong@buaa.edu.cn, sunyq@buaa.edu.cn}\maketitle


 {\begin{center}
\parbox{14.5cm}{\begin{abstract}
 In this paper, a new shift strategy for the implicitly restarted generalized second-order Arnoldi (GSOAR) method is proposed. In implicitly restarted processes, we can get a $k$-step GSOAR decomposition from a $m$-step GSOAR decomposition by performing $p = m-k$ implicit shifted QR iterations. The problem of the implicitly restarted GSOAR is the mismatch between the number of shifts and the dimension of the subspace. There are $2p$ shifts for $p$ QR iterations. We use the shifts to filter out the unwanted information in the current subspace; when more shifts are used, one obtains a better updated subspace. But, if we use more than $p$ shifts, the structure of the GSOAR decomposition will be destroyed. We propose a novel method which can use all $2p$ candidates and preserve the special structure. The new method vastly enhances the overall efficiency of the algorithm. Numerical experiments illustrate the efficiency of every restart process.\vspace{-3mm}
\end{abstract}}\end{center}}



\renewcommand{\baselinestretch}{1.2}
\begin{center} \renewcommand{\arraystretch}{1.5}
{\begin{tabular}{lp{0.8\textwidth}} \hline \scriptsize
{\bf Citation:}\!\!\!\!& \scriptsize Gong F H, Sun Y Q. A new shift strategy for the implicitly restarted generalized second-order Arnoldi method
(in Chinese). Sci Sin Math, 2017, 47: 1{16, doi: 10.1360/012016-22.  You can see the Chinese version online: http://engine.scichina.com/doi/10.1360/012016-22}
 \\
\hline
\end{tabular}}\end{center}

\baselineskip 11pt\parindent=10.8pt  \wuhao
\section{Introduction}

The quadratic eigenvalue problem (QEP)
\begin{equation}
Q(\lambda)x=(\lambda^2 M + \lambda C +K)x=0 \label{01-qep}
\end{equation}
 has various applications in the engineering fields, like dynamic analysis of acoustic systems, fluid mechanics and microelectronic mechanical systems \cite{tisseur,betcke,ruhe}. There are two major classes of numerical methods to solve large QEPs. One is to linearize the QEP into an equivalent generalized eigenvalue problem (GEP) such as
 \begin{equation}
\left [ \begin{array}{cc} -C & -K \\I & 0\end{array} \right ] \left
[\begin{array}{c} \lambda x \\x \end{array} \right ]= \lambda  \left
[ \begin{array}{cc} M & 0 \\0 & I \end{array} \right ] \left [
\begin{array}{c} \lambda x \\x \end{array} \right ]. \label{5xian}
\end{equation}
Currently, the GEP techniques have been quite mature \cite{baidemmel,jia}. But the linearization methods have suffered some disadvantages: losing the spectral property of original QEP and doubling the dimension of the original one. Another class of methods work on the QEP directly. Numerically stable projection methods can keep the spectral properties as well as the structure of original QEP.

The Krylov subspace method plays a significant role in numerical techniques for solving large-scale GEPs. Bai and Su developed the standard Krylov subspace to second-order and proposed the second-order Arnoldi (SOAR) method for large-scale sparse QEPs \cite{baisu}.
They propose a SOAR procedure that computes an orthonormal basis of a second-order Krylov subspace. The SOAR method then projects (\ref{01-qep}) onto this subspace and computes the Ritz pairs to approximate the desired eigenpairs of (\ref{01-qep}). The general convergence property of Rayleigh-Ritz has been demonstrated in \cite{huangt} and the refined Rayleigh-Ritz orthonormal projection method also has been proposed.

On account of the limitation of high computational cost and undue storage requirement, explicit or implicit restarting is usually necessary when we do large-scale matrices calculating. Implicit restarting is a powerful and extensively used technique in eigenvalue problem calculations \cite{sorensen}. Although the SOAR method plays an important role in solving QEPs, we can not directly apply the implicit restarting scheme to the SOAR procedure due to the special structure of the initial vector \cite{meerbergen,otto}. Based on Bai and Su's work, many researchers work specifically to implicitly restart the SOAR procedure. Otto proposed a generalized SOAR (GSOAR) method which substitutes general vector for the original special-structure starting vector. The generalized SOAR procedure can directly apply implicit restarting technique under the hypothesis of no deflation \cite{otto}. An explicit restarting general Krylov subspace method was put forward after making a similar modification by Zhou \cite{zhou}. Huang brought forward the implicit restarting semi orthogonal generalized Arnoldi (SGA) method \cite{huangwq}. At the same time, other researchers are constantly improving the SOAR method in all aspects \cite{jiasun}.

Among these modified restarting SOAR methods, they all encountered the problem that implicit restarting failed to work when shifts did not match up with the dimensions of subspace. From recent researches, it can be seen that the quantity and quality of shifts are of great significance \cite{otto,jiasun}. Jia and Sun explored more properties and features of the GSOAR and gave an efficient method to compute refined Ritz vectors. They proposed a refined GSOAR (RGSOAR) method.
Then they advanced certain exact shifts and refined shifts for respective use within the implicitly restarted GSOAR
and RGSOAR algorithms. The refined shifts are based on the refined Ritz vectors.
They presented an efficient algorithm to compute the exact and refined shift candidates reliably. Unlike the implicitly restarted algorithms for the linear eigenvalue problem, both exact and refined shift candidates are more than the shifts allowed. So the authors showed how to reasonably select the desired shifts among them. In addition, they proposed an effective approach to cure deflation in implicit restarts, so that implicit restarting is useable unconditionally.

Based on Jia and Sun's work, we give a new shift strategy for the implicitly restarted GSOAR and RGSOAR methods which can use all shifts.
In the implicit restarting, we use the unwanted approximate eigenvalues as shifts to filter out the corresponding eigenvector information from the current subspace. A general result is that with more shifts approximating the unwanted eigenvalues, we can generate a better updated subspace.
For the QEPs, one eigenvector can correspond to two different eigenvalues. If the eigenvalue pairs corresponding to the same eigenvector are present in exact or refined shifts, this means some essential shifts will be abandoned. This can be avoided by using all of the shift candidates. Obviously this will destroy the original structure. We use a subtle transformation making the damaged structure back to the Hessenberg form. Due to all shifts, the unwanted information is excluded more thoroughly and we can get a better subspace in each restarting. Because the implicit restarting in \cite{jiasun} is useable unconditionally, our algorithm is useable inherently.

The rest of this paper is arranged as follows. In section 2, we introduce the SOAR method and GSOAR method, implicit restarting,  selection of exact and refined shifts. In section 3, we present the new strategy: using all shifts for projection method to enhance the efficiency of the method. In section 4, several numerical experiments are presented to illustrate the efficiency of the new implicitly restarted GSOAR method.
\\
\\
\section{The implicitly restarted GSOAR and RGSOAR methods }

The SOAR method was proposed to solve the large-scale sparse QEPs based on the second-order Krylov subspace by Bai and Su \cite{baisu}.
\begin{definition}\quad
Let matrices $A$, $B\in \mathcal{C}^{n\times n}$, non-zero vector $u\in \mathcal{C}^n$, and define the sequence $r_0,r_1,r_2,\cdots,r_{n-1}$ based on $A$, $B$ and $u$, where
$$
\begin{array}{ccl}
r_0&=&u\\
r_1&=&Au\\
r_j&=&Ar_{j-1}+Br_{j-2} \hspace{0.5cm} \mbox{for } j\geq 2.
\end{array}
$$
The space \begin{equation}\label{4}
\mathcal{G}_{n}(A,B;u)=span\{r_0,r_1,r_2,\cdots,r_{n-1}\}
\end{equation} 
is called a $n$th second-order Krylov subspace.
\end{definition}

If we apply standard Krylov subspace technique to \eqref{5xian}, then the corresponding Krylov subspace is:
\begin{equation}\label{5}
\kappa_n(H,\nu)=span\{\nu,H\nu,H^{2}\nu,\cdots,H^{n-1}\nu\},
\end{equation} 
\noindent with $\nu=
\begin{bmatrix}
u \\
0
\end{bmatrix}
$.
After some simple calculations, we note that there are some associations between the standard Krylov vectors $H^{j}\nu$ of length $2n$ in \eqref{5}
\begin{equation}\label{6}
H^{j}\nu=
\begin{bmatrix}
r_{j}  \\
r_{j-1}
\end{bmatrix},j\geq1,
\end{equation}
and the vectors of second-order Krylov subspace
\begin{equation}\label{7}
r_j=Ar_{j-1}+Br_{j-2}.
\end{equation} 

The subspace \eqref{5} based on matrix $H$ and starting vector $v$ can be generalized by vector sequence $r_0,r_1,r_2,\cdots,r_{n-1}$. In the meantime, \eqref{6} indicates that the second-order Krylov subspace $\mathcal{G}_{j}(A,B;u)$ can be used as the projection subspace of QEPs while we only need to generate an orthonormal basis ${q_j}_{j=1}^n$:
\begin{equation}\label{8}
\mathcal{G}_{j}(A,B;u)=span\{q_1,q_2,\cdots,q_{j}\},j\geq1.
\end{equation} 

Bai and Su presented a procedure for generating an orthonormal basis of the second-order Krylov subspace $\mathcal{G}_{j}(A,B;u)$.

\begin{alg} {\bf SOAR procedure with deflation remedy}
\begin{enumerate}
\item[1:] $q_1=u_1/\parallel u_1\parallel_2$
\item[2:] $p_1=0$
\item[3:]  for $j=1:m$
\item[4:] \hspace{0.5cm}$r=Aq_j+Bp_j$
\item[5:] \hspace{0.5cm} $s=q_j$
\item[6:] \hspace{0.5cm}for $i=1:j$
\item[7:] \hspace{1cm} $t_{ij}=q_i^{T}r$
\item[8:] \hspace{1cm} $r=r-q_it_{ij}$
\item[9:] \hspace{1cm}$s=s-p_it_{ij}$
\item[10:] \hspace{0.5cm}end for
\item[11:] \hspace{0.5cm}$t_{j+1,j}=\parallel r\parallel_2$
\item[12:] \hspace{0.5cm}if $t_{j+1,j}=0$,
\item[13:] \hspace{1cm}if $s\in span\{p_i\mid i:q_i=0,i\leq i\leq j\}$
\item[14:] \hspace{1.5cm} break
\item[15:] \hspace{1cm}else deflation
\item[16:] \hspace{1.5cm} $t_{j+1,j}=1$
\item[17:] \hspace{1.5cm} $q_{j+1}=0$
\item[18:] \hspace{1.5cm} $p_{j+1}=s$
\item[19:] \hspace{1cm} end if
\item[20:] \hspace{0.5cm}else
\item[21:] \hspace{1cm}$p_{j+1}=s/t_{j+1,j}$
\item[22:] \hspace{1cm} $q_{j+1}=r/t_{j+1,j}$
\item[23:] \hspace{0.5cm} end if
\item[24:] end for
\end{enumerate}
\end{alg}

At step $j$, if $r_i,\, i=0,1, \ldots, j$ are linearly dependent but $[r_i^T,r_{i-1}^T]^T,\,i=0, \ldots, j$
with $r_{-1}=0$ are not, we call this situation deflation; if both $\{r_i\}$ and ~$\{[r_i^T,r_{i-1}^T]^T\}$ are
linearly dependent at step $j$, we call this situation breakdown. If $t_{j+1,j}=0$ at a certain step $j$ $(j<k)$, that means deflation occurs. In this case the second-order Krylov subspace $\mathcal{G}_{j}(A,B;u)$ does not contain any exact eigenvector of \eqref{01-qep} \cite{jiasun}. Therefore, the measure to remedy deflation is to reset $t_{j+1,j}$ a non-zero value, here simply set it to be one, let $q_{j+1}=0$, and continue the algorithm. The non-zero vectors of sequence $\{q_j\}$  still maintain orthonormal and span the subspace $\mathcal{G}_{j}(A,B;u)$.

For the quadratic eigenvalue problem \eqref{01-qep} of $n$ dimension , it corresponds to a second-order Krylov subspace $\mathcal{G}_{m}(A,B;u)$, with $m\leq n$, and $A=-M^{-1}C$, $B=-M^{-1}K$. Let $Q_m=[q_1,q_2,\cdots,q_m]\in \mathcal{C}^{n\times m}$, $P_m=[p_1,p_2,\cdots,p_m]\in \mathcal{C}^{n\times m}$ and unreduced upper Hessenberg matrix $\hat{T}_{m}=[{t_{ij}} ]=\begin{bmatrix}
   {{T}_{m}}  \\
   {{t}_{m+1,m}}e_{m}^{T}  \\
\end{bmatrix} \in \mathcal{C}^{(m+1)\times m}$.
If Algorithm 1 does not stop before $m$, then we have
\begin{equation}\label{9}
span\{Q_m\}=\mathcal{G}_{m}(A,B;u),
\end{equation} 
and define
\begin{equation}\label{10}
H \begin{bmatrix}
   {{Q}_{m}}  \\
   {{P}_{m}}  \\
\end{bmatrix}= \begin{bmatrix}
   {{Q}_{m+1}}  \\
   {{P}_{m+1}}  \\
\end{bmatrix} \hat{T}_{m}
\end{equation}
as $m$-step SOAR decomposition, where
$Q_{m+1}=[Q_m,q_{m+1}]$, $P_{m+1}=[P_m,p_{m+1}]$, $H=\begin{bmatrix}
A & B  \\
I & 0
\end{bmatrix}$.

Rayleigh-Ritz method can be adapted in \eqref{01-qep} to seek the approximate eigenpairs $(\theta,y)$ satisfying the Galerkin restriction:
\begin{equation}\label{11}
(\theta ^{2}M+\theta C+K)y\bot \mathcal{G}_m(A,B;u),
\end{equation}
\noindent where $y\in \mathcal{G}_m(A,B;u)$, $\theta\in \mathcal{C}$.  The columns of $n\times n$ matrix $Q_m$ generated by Algorithm 1 are an orthonormal basis of $\mathcal{G}_m(A,B;u)$. Let $\theta$, $g$  satisfy the small-scale QEP:
\begin{equation}\label{12}
(\theta^2M_m+\theta C_m+K_m)g=0,
\end{equation}
\noindent where $M_m=Q_m^{T}MQ_m$, $C_m=Q_m^{T}CQ_m$, $K_m=Q_m^{T}KQ_m$, then we can obtain the approximate eigenpairs $\theta$ and
$y=Q_{m}g$ of \eqref{01-qep} .
The small-scale QEP \eqref{12} is produced by projecting the original QEP \eqref{01-qep} onto $\mathcal{G}_m(A,B;u)$ explicitly. This  method is named second-order Arnoldi (SOAR) method since it is an Arnoldi-like procedure. From the theory of Arnoldi method, we can find that the SOAR method has the advantages of fast convergence rate and simultaneous convergence of a group of eigenvalues. Moreover, it can still keep the special structure of original QEP after projection.

A disadvantage of this method is that implicit restarting scheme cannot be directly adapted. The new starting vector $p_1$ will no longer be zero after truncating. In order to satisfy the implicitly restarted conditions, researchers suggested the GSOAR method which replaces $p_1=0$ by $p_1=u_2/\parallel u_2\parallel_2$ and meets the requirement that starting vector $p_1=u_2/\parallel u_2\parallel_2$ is non-zero, where $u_2$ is a non-zero vector. Combining the deflation remedy SOAR procedure with GSOAR procedure, we get the following algorithm:

\begin{alg} {\bf GSOAR procedure with deflation remedy}\label{GSOAR}
\begin{enumerate}
  \item[1:] $q_1=u_1/\parallel u_1\parallel_2$, $p_1=u_2/\parallel u_2\parallel_2$.
  \item[2:]for $j=1,2,\ldots,k$ do
  \item[3:] \hspace{0.5cm} $r=Aq_j+Bp_j$
  \item[4:]\hspace{0.5cm}  $s=q_j $
  \item[5:]\hspace{0.5cm}  for ~$i=1,2,\ldots,j$ do
  \item[6:]\hspace{1cm}     $t_{ij}=q_i^*r $
  \item[7:]\hspace{1cm}     $r=r-t_{ij}q_i$
  \item[8:]\hspace{1cm}     $s=s-t_{ij}p_i$
  \item[9:]\hspace{0.5cm}  end for
  \item[10:]\hspace{0.5cm}   $t_{j+1j}=\|r\|$
  \item[11:]\hspace{0.5cm}  if  ~$t_{j+1j}=0$
  \item[12:]\hspace{1cm}    if $s\in span\{p_i|i:
                             q_i=0,1\leq i\leq j  \}$
  \item[13:]\hspace{1.5cm}   break
  \item[14:]\hspace{1cm}   else   {\em deflation}
  \item[15:]\hspace{1.5cm} reset $t_{j+1j}=1$
  \item[16:]\hspace{1.5cm}  $q_{j+1}=0$
  \item[17:]\hspace{1.5cm}   $p_{j+1}=s$
  \item[18:]\hspace{1cm}     end if
  \item[19:]\hspace{0.5cm}   else
  \item[20:]\hspace{1cm }    $q_{j+1}=r/t_{j+1j}$
  \item[21:]\hspace{1cm }    $p_{j+1}=s/t_{j+1j}$
  \item[22:]\hspace{0.5cm }   end if
  \item[23:]   end for
\end{enumerate}
\end{alg}

If Algorithm \ref{GSOAR} does not stop before $m$, it gives the following generalized second-order
Krylov sequence and subspace; see \cite{jiasun}.

Let $A$ and $B$ be $n\times n$ matrices and for vectors $u_1, u_2\in
\mathcal{C}^n$, and define
$$
\begin{array}{ccl}
r_0&=&u_1,\\
r_1&=&Ar_0+Bu_2,\\
r_j&=&Ar_{j-1}+Br_{j-2} \hspace{0.5cm} \mbox{for } j\geq 2.
\end{array}
$$
Then $r_0, r_1, r_2, \ldots, r_{m-1}$
is called a generalized second-order Krylov sequence
based on $A,B$ and $u_1,u_2$,
and
$
\mathcal{G}_m(A,B;u_1,u_2)={span}\{r_0, r_1, r_2, \ldots, r_{m-1}\}
$
the $m$-th generalized second-order Krylov subspace.

Let $Q_m=[q_1,q_2,\cdots,q_m]\in \mathcal{C}^{n\times m}$, $P_m=[p_1,p_2,\cdots,p_m]\in \mathcal{C}^{n\times m}$ and
$\hat{T}_{m}=[{t_{ij}} ]=\begin{bmatrix}
   {{T}_{m}}  \\
   {{t}_{m+1,m}}e_{m}^{T}  \\
\end{bmatrix} \in \mathcal{C}^{(m+1)\times m}$.
Then $Q_m$ is an orthonormal basis of  the $m$-th generalized second-order Krylov subspace $\mathcal{G}_{m}(A,B;u_1,u_2)$
and we get the $m$-step GSOAR decomposition:
\begin{equation}\label{10}
H \begin{bmatrix}
   {{Q}_{m}}  \\
   {{P}_{m}}  \\
\end{bmatrix}= \begin{bmatrix}
   {{Q}_{m+1}}  \\
   {{P}_{m+1}}  \\
\end{bmatrix} \hat{T}_{m},
\end{equation}
where
$Q_{m+1}=[Q_m,q_{m+1}]$, $P_{m+1}=[P_m,p_{m+1}]$, $H=\begin{bmatrix}
A & B  \\
I & 0
\end{bmatrix}$.

We can project the original QEP \eqref{01-qep} onto the generalized second-order Krylov subspace $\mathcal{G}_{m}(A,B;u_1,u_2)$ explicitly.
We still get the  small-scale QEP:
\begin{equation}\label{12new}
(\theta^2M_m+\theta C_m+K_m)g=0,
\end{equation}
\noindent where $M_m=Q_m^{T}MQ_m$, $C_m=Q_m^{T}CQ_m$, $K_m=Q_m^{T}KQ_m$.
We can get the Ritz pairs $\theta$ and $y=Q_{m}g$ of \eqref{01-qep} by computing the eigenpairs of \eqref{12new}.
Then we get the GSOAR method.
Suppose that we have computed the Ritz values $\theta$ by the
GSOAR method. For each $\theta$, we seek a unit length vector
$\tilde{u} \in \mathcal{G}_m(A, B;u_1,u_2)$
satisfying the optimal requirement
\begin{equation}
\|(\theta^2 M + \theta C+ K) {\tilde{u}}\| = \arg \min_{
 \mbox{\scriptsize $\begin{array}{c} u \in \mathcal{G}_m(A,
B;u_1,u_2) \\ \|u\|=1
 \end{array}$}}
\|(\theta^2 M + \theta C+ K )u\| \label{th2}
\end{equation}
and use it as an approximate eigenvector, called the refined Ritz
vector. The pairs $(\theta,\tilde u)$ are also called the refined Rayleigh--Ritz
approximations and the method is called refined generalized SOAR (RGSOAR) method.
The refined vector can be computed by seeking a unit length vector
$\tilde{z} \in {\cal C}^m$ such that
$\tilde{u}=Q_m\tilde z$
with
\begin{equation}
\tilde z=\arg\min_{ \mbox{\scriptsize $\begin{array}{c} z \in {\cal C}^m \\
\|z\|=1 \end{array}$}} \|(\theta^2 M + \theta C+ K)Q_mz\|,
\label{th3}
\end{equation}
the right singular vector of the matrix
$(\theta^2 MQ_m+\theta CQ_m+KQ_m)$ associated with its smallest
singular value $\sigma_{\min} (\theta^2 MQ_m+\theta CQ_m+KQ_m)$.
It was shown in \cite{jiasun} that the singular vectors can be computed with moderate cost.

\subsection{Implicit restarting}

For the SOAR method and GSOAR method, in order to guarantee the convergence of the approximate eigenpairs, $m$ should be as large as possible theoretically. Nevertheless, the storage and computing requirements become massive as the dimension $m$ of the subspace increases. So restarting is necessary for a practical method. If no deflation occurs, it is direct to adapt the implicit restarting scheme \cite{sorensen} to the GSOAR procedure. Given $p$ shifts $\mu_1,\mu_2,\ldots,\mu_p$, we can apply $p$ implicit shifted QR algorithm to matrix $T_m$ in \eqref{10}
\begin{equation}\label{14}
(T_m-\mu_{1}I)\cdots(T_m-\mu_{p}I)=V_m R,
\end{equation}
\noindent where $V_m$ is a $m\times m$ orthogonal matrix and $R$ is a $m\times m$ upper triangular matrix. By relation \eqref{10}, we have
\begin{equation}\label{15}
H \begin{bmatrix}
   {{Q}_{m}}  \\
   {{P}_{m}}  \\
\end{bmatrix}V_m= \begin{bmatrix}
   {{Q}_{m}}  \\
   {{P}_{m}}  \\
\end{bmatrix} V_m (V_m^T T_m V_m)+t_{m+1,m}\begin{bmatrix}
{{q}_{m+1}}  \\
   {{p}_{m+1}}  \\
   \end{bmatrix}e_m^T V_m.
\end{equation}

\noindent With simple computation, it is noted that $V_m^T T_m V_m$ is a Hessenberg matrix. If we take
$\begin{bmatrix}
   {{\tilde{Q}}_{m}}  \\
   {{\tilde{P}}_{m}}  \\
\end{bmatrix}=\begin{bmatrix}
   {{Q}_{m}}  \\
   {{P}_{m}}  \\
\end{bmatrix} V_m$
and take $V_m^T T_m V_m=\tilde{T}_m$, then we have
\begin{equation}\label{16}
H\begin{bmatrix}
   {{{\tilde{Q}}}_{m}}  \\
   {{{\tilde{P}}}_{m}}  \\
\end{bmatrix} =\begin{bmatrix}
   {{{\tilde{Q}}}_{m}}  \\
   {{{\tilde{P}}}_{m}}  \\
\end{bmatrix} {{\tilde{T}}_{m}}+{{t}_{m+1,m}} \begin{bmatrix}
   {{q}_{m+1}}  \\
   {{p}_{m+1}}  \\
\end{bmatrix} b_{m}^{T},
\end{equation}
\noindent where $b_{m}^{T}=( \begin{matrix}
   0 & \cdots  & 0 & \begin{matrix}
   \begin{matrix}
   {{b}_{m-p}} & \cdots   \\
\end{matrix}b_{m-1}^{{}} & {{b}_{m}}  \\
\end{matrix}  \\
\end{matrix}) $
which contains at least $m-p-1$ zeros.

Given a $m$-step GSOAR decomposition \eqref{16}, it can be truncated as a $k$-step GSOAR decomposition while $k\leq m-p$
\begin{equation}\label{17}
H\begin{bmatrix}
   {{{\tilde{Q}}}_{k}}  \\
   {{{\tilde{P}}}_{k}}  \\
\end{bmatrix} =\begin{bmatrix}
   {{{\tilde{Q}}}_{k}}  \\
   {{{\tilde{P}}}_{k}}  \\
\end{bmatrix} {{\tilde{T}}_{k}}+{{t}_{k+1,k}} \begin{bmatrix}
   {{\tilde{q}}_{k+1}}  \\
   {{\tilde{p}}_{k+1}}  \\
\end{bmatrix} e_{k}^{T}.
\end{equation}
\noindent
Meanwhile, we get a $k$-th restarted subspace $G_{k}(A, B;u_1, u_2)$ spanned by an orthonormal vectors sequence $\{\tilde{q_1},\tilde{q_2},\cdots,\tilde{q_k}\}$.

Because the implicitly restarted GSOAR method has a strict requirement for the dimension of restarted subspace, we can truncate \eqref{16} into a new $k$-th restarted subspace if and if only the subspace dimension $k\leq m-p$. Specifically, for \eqref{17}, when $k$ does not exceed $m-p$, the last term in the right looks like $e_k^T$ after truncating and \eqref{17} will have the same form as \eqref{10}. The $k$-step GSOAR decomposition can be expanded to $m$-step GSOAR decomposition again via Algorithm \ref{GSOAR}. Then looking for new Ritz pairs, computing their deviation and the process repeats.

Equation \eqref{17} is equivalent to a $k$-step  GSOAR  procedure mathematically. If deflations occur at steps $m_1,m_2,\ldots,m_j\leq k$, then the corresponding $j$ columns $\tilde{q}_{m_j}$ of $\tilde{Q}_{k}$ are zeros.
Jia and Sun gave a method to cure this problem, so the implicit restarting is unconditional\ref{jiasun}.

\subsection{Exact shifts and refined shifts}

We consider the computation of shifts for GSOAR and RGSOAR methods in this subsection.
Two kinds of shifts were given in\cite{jiasun}: the exact shifts and refined shifts.
We first show how to compute the refined shifts as example.

The RGSOAR method computes the refined Ritz vectors
$\tilde{u}_i$, which can be much more accurate
than the Ritz vectors $y_i$.There were certain refined shifts
for the refined Arnoldi method and the refined harmonic Arnoldi method for the
linear eigenvalue problem\cite{jia,jia02}. In the same spirit, the
refined shifts were proposed for the RGSOAR method.

Let $span\{\hat{u}_1,\ldots,\hat{u}_p\}$ be the orthogonal complement of the refined Ritz vectors $\tilde{u}_i,\,i=1,2,\ldots,k$ with respect to $\mathcal{G}_k(A,B;q_1,p_1)$.
Project QEP (\ref{01-qep}) onto the $span\{\hat{u}_1,\ldots,\hat{u}_p\}$. Then we obtain a $p=m-k$-dimensional projected QEP and compute its
$2p$ eigenvalues. As a result, the $2p$  eigenvalues are the unwanted approximate eigenvalues.
We call them refined shift candidates.
We select $p=m-k$ refined shift candidates as shifts, called the refined shifts,
for use within the implicitly restarted RGSOAR algorithm.

For GSOAR method, we can take the $2(m-k)$ unwanted Ritz values as the shift candidates directly.
But it is just possible that one wanted Ritz value and one unwanted Ritz value correspond to the same Ritz vector.
When we take the unwanted Ritz value as shift, that will filter out the wanted Ritz vector by restarting.
To avoid this situation, we can compute the exact shift candidates in the orthogonal complement of the Ritz vectors $y_i,\,i=1,2,\ldots,k$ with respect to $\mathcal{G}_k(A,B;q_1,p_1)$. We select $p=m-k$ exact shift candidates as shifts, called the exact shifts, for GSOAR method.

Both methods have $2p$ shift candidates, but if we use more than $p$ shifts,  equation (\ref{16})  can not be truncated at the first $k$ columns. So we can select the $p$ shifts farest from the target as shifts.

\section{Shift strategy for implicitly restarted GSOAR and RGSOAR methods}

In both methods, all refined vectors $\tilde{u}_1,\ldots,\tilde{u}_k$ and Ritz vectors ${y}_1,\ldots,{y}_k$ are approximating the desired eigenvectors $x_1,\ldots,x_k$. The exact and refined shift candidates are
approximations to some of the unwanted eigenvalues of QEP (\ref{01-qep}) because the information on $x_{k+1},\ldots,x_m$ has been
removed from $\mathcal{G}_m(A,B;q_1,p_1)$. After restarting in (\ref{17}), RGSOAR obtains  $span\{\tilde{Q}_k\}=span\{\tilde{u}_1,\ldots,\tilde{u}_k$\} and GSOAR obtains  $span\{\tilde{Q'}_k\}=span\{y_1,\ldots,y_k$\}. Both methods achieve this goal through removing the unwanted information by shifts. As mentioned above, if two different shifts corresponding to same unwanted eigenvalue, at least one unwanted eigenvector can not be removed from the subspace $\mathcal{G}_m(A,B;q_1,p_1)$. The result of restarting can not be guaranteed. In order to achieve the restarting goal, all shift candidates must be used.
The relationship between the number of shifts and the number of zero elements of $b_{m}^{T}=( \begin{matrix}
   0 & \cdots  & 0 & \begin{matrix}
   \begin{matrix}
   {{b}_{m-p}} & \cdots   \\
\end{matrix}b_{m-1}^{{}} & {{b}_{m}}  \\
\end{matrix}  \\
\end{matrix}) $
in \eqref{16} is shown in Figure 1. If  the number of shifts $p$ is more than $m-k$, the first $k$ columns of  \eqref{16} will  not be a GSOAR decomposition. Because the last term in the right side is not only with $e_k^T$ but also with $e_{k-1}^T$, $e_{k-2}^T$ until $e_{m-p}^T$. So the implicit restarting technique can not be applied directly.
Our main work is to propose an implicitly restarted generalized second-order Arnoldi strategy which can use all shift candidates.

It was shown in \cite{stewart}, that Krylov decomposition and Arnoldi decomposition are equivalent. If we can take Figure 1 (d) as a generalized Krylov decomposition, then it can be transformed to a GSOAR decomposition by orthogonal transformation. The scheme is written as the following result.

\begin{figure}[!h]
\begin{center}
\subfigure[GSOAR decomposition]{
\resizebox*{5cm}{!}{\includegraphics{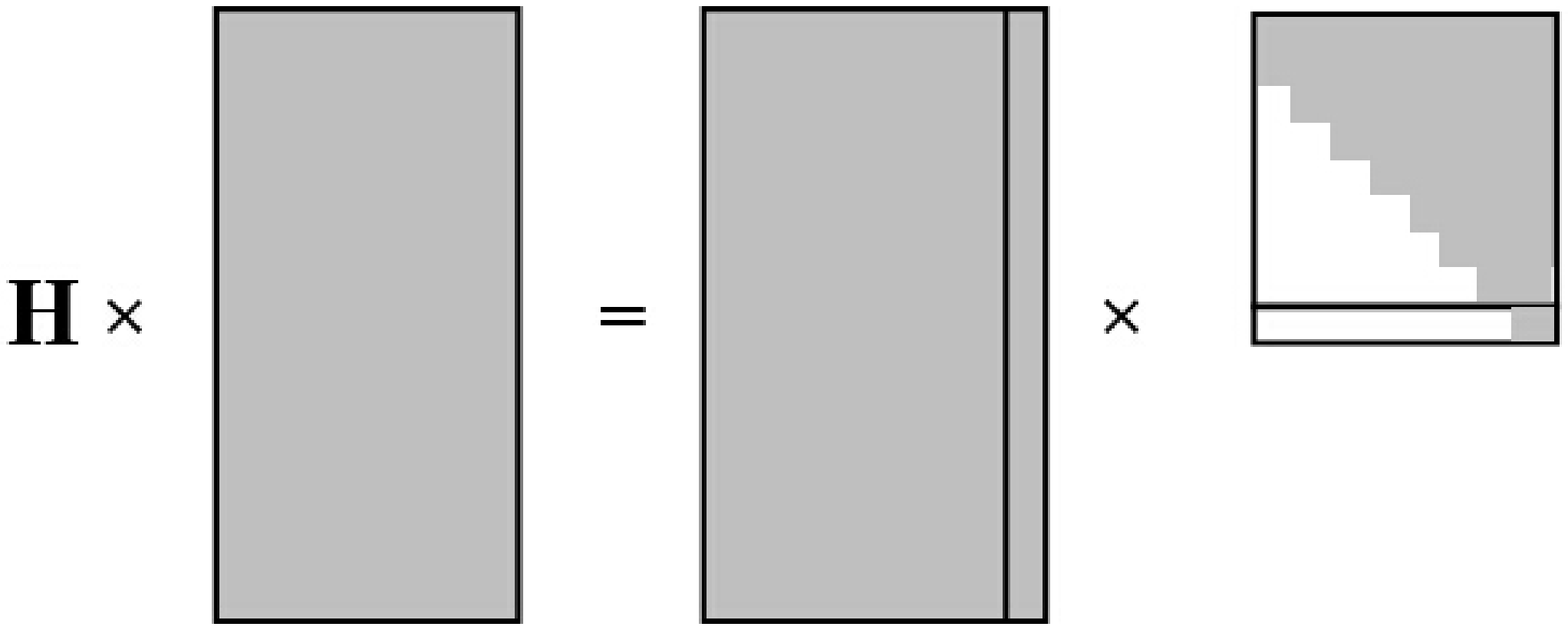}}}\hspace{5pt}
\subfigure[$p=m-k$]{
\resizebox*{5cm}{!}{\includegraphics{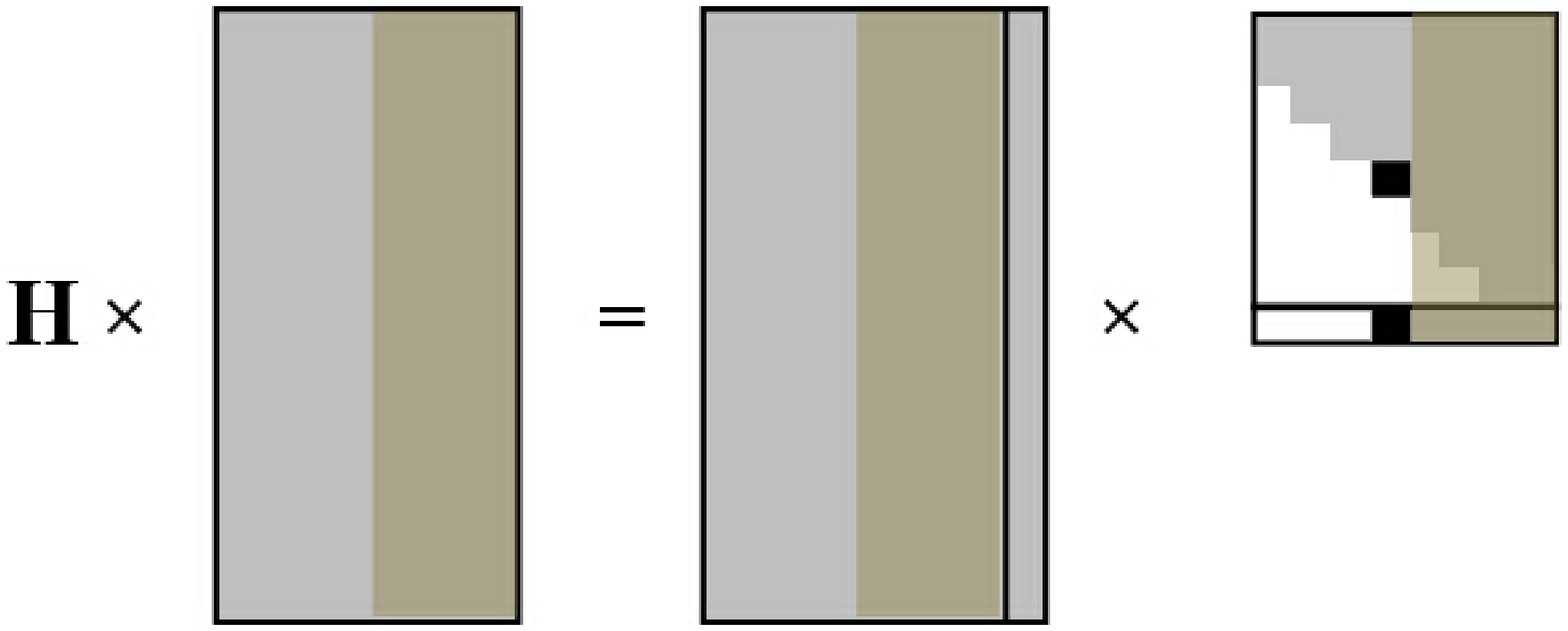}} }\\
\subfigure[$p>m-k$]{
\resizebox*{5cm}{!}{\includegraphics{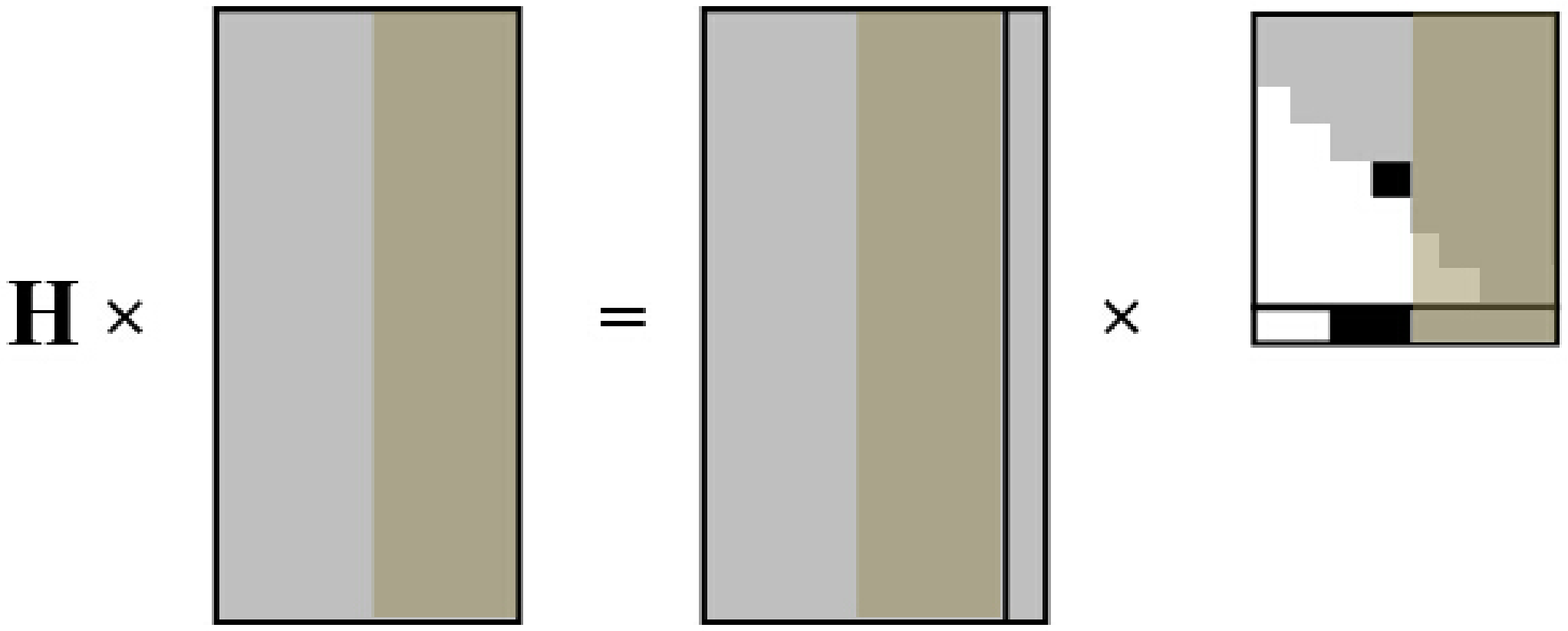}}}\hspace{5pt}
\subfigure[$p>m$]{
\resizebox*{5cm}{!}{\includegraphics{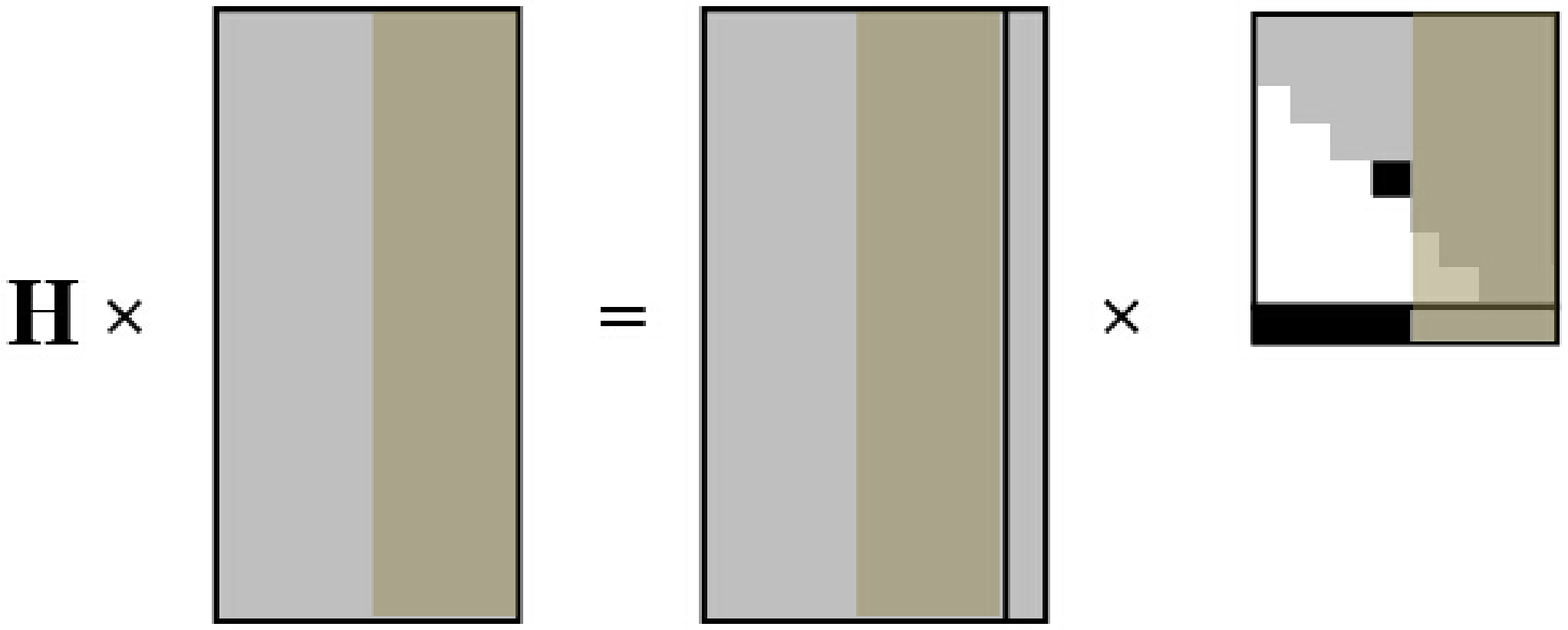}}} \\
\caption{Implicitly restarted GSOAR method schematic}
\label{pic2}
\end{center}
\end{figure}

\begin{theorem}\quad
For equation \eqref{16}, there exists an orthonormal matrix $W$, that transforms $b_m^T$ into $\alpha e_m^T$, where $\alpha$ is a non-zero constant, and $W^T \tilde{T}_m W$ is a Hessenberg matrix.
\end{theorem}

\begin{proof}\quad
Let $W_1$ be a Householder matrix, such that
\begin{equation}\label{22or33}
W_1b_{m}=\alpha e_{m},
\end{equation} 
\noindent with $\alpha$ as a non-zero constant. Transpose of above equation,
\begin{equation}\label{23or34}
b_{m}^{T}W_1^{T}=\alpha e_{m}^{T},
\end{equation} 
with $W_1=W_1^T$, we have
\begin{equation}\label{24or35}
b_{m}^{T}W_1=\alpha e_{m}^{T}.
\end{equation} 
Then multiply matrix $W_1$ on the right side of equation \eqref{16}
\begin{equation}\label{25or36}
H
\begin{bmatrix}
\widetilde{Q}_{m}  \\
\widetilde{P}_{m}
\end{bmatrix}
W_1=
\begin{bmatrix}
\widetilde{Q}_{m}  \\
\widetilde{P}_{m}
\end{bmatrix}
W_1(W_1^{T}\widetilde{T}_{m}W_1)+t_{m+1,m}
\begin{bmatrix}
q_{m+1} \\
p_{m+1}
\end{bmatrix}
b_{m}^{T}W_1 .
\end{equation} 
Take $\begin{bmatrix}
\widetilde{Q}_{m}  \\
\widetilde{P}_{m}
\end{bmatrix}W$
as a new $\begin{bmatrix}
Q_m \\
P_m
\end{bmatrix}$
and denote $W_1^{T}\widetilde{T}_{m}W_1$ by new matrix $B_m$. Note that $B_m$ loses the Hessenberg form. Combining the above two equations, we have
\begin{equation}\label{26or37}
H
\begin{bmatrix}
Q_{m}  \\
P_{m}
\end{bmatrix}
=
\begin{bmatrix}
Q_{m}  \\
P_{m}
\end{bmatrix}
B_m+b_{m+1,m}
\begin{bmatrix}
q_{m+1}  \\
p_{m+1}
\end{bmatrix}
e_{m}^{T} ,
\end{equation} 
where $ b_{m+1,m} = t_{m+1,m}\alpha $. At this moment, if we want to do implicit restarting, we must transform $B_m$ into Hessenberg again and maintain the last term $e_{m}^{T}$.

Householder transformation and orthogonal transformation techniques are adopted once again to transform $B_m$ into upper Hessenberg matrix. To prevent the next transformation damaging the previous result, the Householder transformations are performed on $B_m$ from the bottom row to the top row.

When $l=1$, denote $B_m$ by matrix block form:
\begin{equation}\label{27or38}
B_m=
 \begin{bmatrix}
  b_{11} & b_{12} & \cdots & b_{1m} \\
  b_{21} & b_{22} & \cdots & b_{2m} \\
  \vdots  & \vdots  & \ddots & \vdots  \\
  b_{m1} & b_{m2} & \cdots & b_{mm}
 \end{bmatrix}
 =
 \begin{bmatrix}
B_{11} & B_{12} \\
c_1 & b_{mm}
\end{bmatrix},
\end{equation} 
\noindent where $
c_1=\begin{pmatrix}
b_{m,1}  & \cdots & b_{m,m-1}
\end{pmatrix}$, $B_{11}= B_m(1:m-1,1:m-1)$ is a $(m-1)\times (m-1)$ matrix, $B_{12}=B_m(1:m-1,m)$ is a column vector of length $m-1$. Let $R_1$ be a $(m-1)\times (m-1)$ Householder matrix such that
\begin{equation}\label{28or39}
c_1R_1=\alpha_1 e_{m-1}^T,
\end{equation} 
where $\alpha_ 1$ is a non-zero constant. Let $U_1=\begin{bmatrix}
R_1 & 0\\
0 & 1
\end{bmatrix}$, we perform orthogonal transformation on $B_m$,
\begin{equation}\label{29or310}
B_{m}^{(1)}=U_{1}^{T} B_m U_{1},
\end{equation}
then $B_m^{(1)}$  has the following structure
\begin{equation}
 U_{1}^{T}
 \begin{bmatrix}
B_{11} & B_{12} \\
c_1 & b_{mm}
\end{bmatrix}
U_1=
\arraycolsep5pt
\left[
\begin{array}{@{\,}ccc|c@{\,}}
&B_{11}^{(1)}&&B_{12}^{(1)}\\
\hline
&c_{2}&&\\
0 &\cdots &0& B_{22}^{(1)}\\
\end{array}
\right],
\end{equation} 
where $B_{11}^{(1)}$ is a $(m-2)\times (m-2)$ matrix, $B_{12}^{(1)}$ is a $(m-2)\times 2$ matrix, $B_{22}^{(1)}$ is a $2\times 2$ matrix and $ c_2=\begin{pmatrix}
b_{m-1,1}  & \cdots & b_{m-1,m-2}
\end{pmatrix}$.

When $l=2$, we can build a $(m-2)\times (m-2)$  Householder matrix $R_2$, such that
\begin{equation}\label{30or311}
c_2R_2=\alpha_2 e_{m-2}^T ,
\end{equation} 
where $\alpha_ 2$ is a non-zero constant and use
$U_2=\begin{bmatrix}
R_2 & 0\\ 0 & I_{2\times 2}
\end{bmatrix}$ to transform  $B_{m}^{(1)}$ into
\begin{equation}\label{31or312}
B_m^{(2)}=U_{2}^{T} B_m^{(1)} U_{2}=
\arraycolsep5pt
\left[
\begin{array}{@{\,}ccc|c@{\,}}
&\tilde{B}_{11}^{(1)}&&\tilde{B}_{12}^{(1)}\\
\hline
0&\cdots&\alpha_{2}&\\
0 &\cdots &0& B_{22}^{(1)}\\
\end{array}
\right]
=
\arraycolsep5pt
\left[
\begin{array}{@{\,}ccc|c@{\,}}
&B_{11}^{(2)}&&B_{12}^{(2)}\\
\hline
&c_{3}& &\\
0 &\cdots & 0& B_{22}^{(2)}\\
0 &\cdots &0&
\end{array}\right],
\end{equation} 
where
\begin{equation}
B_{22}^{(2)} =
\left[
\begin{array}{@{\,}c|c@{\,}}
\ast &
 \ast, \ast  \\
\hline
\begin{array}{c}
\alpha_2 \\ 0 \end{array} & B_{22}^{(1)}\\
\end{array}\right]
\end{equation}
is a $3\times 3$ Hessenberg matrix.

Suppose $B_{m}^{(k-1)}$ is the following matrix block form:
\begin{equation}\label{32or313}
B_m^{(k-1)}=
\arraycolsep5pt
\left[
\begin{array}{@{\,}ccc|c@{\,}}
&B_{11}^{(k-1)}&&B_{12}^{(k-1)}\\
\hline
&c_{k}&&\\
0 &\cdots &0&\\
\vdots&\cdots&\vdots&B_{22}^{(k-1)} \\
0 &\cdots &0&
\end{array}
\right],
\end{equation} 
\noindent where $B_{11}^{(k-1)}$ is a $(m-k)\times (m-k)$ matrix, $B_{12}^{(k-1)}$ is a $(m-k)\times k$ matrix, $B_{22}^{(k-1)}$ is a $k\times k$  Hessenberg  matrix and $
c_k=\begin{pmatrix}
b_{m-k+1,1}  & \cdots & b_{m-k+1,m-k}
\end{pmatrix}$. We can build a $(m-k)\times (m-k)$ Householder matrix $R_k$, such that
\begin{equation}\label{33or314}
c_kR_k=\alpha_k e_{m-k}^T,
\end{equation} 
where $\alpha_ k$ is a non-zero constant. Let
$U_k=\begin{bmatrix}
R_k & 0\\
0 & I_{k\times k}
\end{bmatrix}$ , we obtain $B_{m}^{(k)}$:
\begin{equation}\label{34or315}
B_m^{(k)}=U_k^TB_m^{(k-1)}U_k=
\arraycolsep5pt
\left[
\begin{array}{@{\,}ccc|c@{\,}}
&\tilde{B}_{11}^{(k-1)}&&\tilde{B}_{12}^{(k-1)}\\
\hline
0&\cdots&\alpha_k&\\
\vdots&\cdots&\vdots&B_{22}^{(k-1)} \\
0 &\cdots &0& \\
\end{array}
\right]
= \arraycolsep5pt
\left[
\begin{array}{@{\,}ccc|c@{\,}}
&B_{11}^{(k)}&&B_{12}^{(k)}\\
\hline
&c_{k+1}&&\\
0 &\cdots &0&\\
\vdots&\cdots&\vdots&B_{22}^{(k)} \\
0 &\cdots &0& \\
\end{array}
\right],
\end{equation} 
with $B_{11}^{(k)}$ is a $(m-k-1)\times (m-k-1)$ matrix, $B_{12}^{(k)}$ is a $(m-k-1)\times (k+1)$ matrix, $B_{22}^{(k)}$ is a $ (k+1)\times (k+1)$ Hessenberg matrix and $
c_k=\begin{pmatrix}
b_{m-k,1}  & \cdots & b_{m-k,m-k-1}
\end{pmatrix}$.

Obviously, when $l = m-2$, $B^{(m-2)}_m$ is an upper Hessenberg matrix.

Let $U=U_1U_2\cdots U_{m-2}$, obviously $U$ is an orthonormal matrix as $U_1,U_1,\cdots ,U_{m-2}$ are orthonormal matrices. Multiply $U$ on equation \eqref{26or37} from right:
\begin{equation}\label{37or318}
H
\begin{bmatrix}
Q_{m}  \\
P_{m}
\end{bmatrix}
U =  \begin{bmatrix}
Q_{m}  \\  P_{m}
\end{bmatrix}
U(U^{T}B_{m}U)+b_{m+1,m}
\begin{bmatrix}
q_{m+1} \\ p_{m+1}
\end{bmatrix}
e_{m}^{T}U .
\end{equation} 
Because $ U_1,U_1,\cdots ,U_{m-2}$ have special structure, the last row of $U$ must be $e_m^T$. Therefore, we have
\begin{equation}\label{38or319}
e_m^{T}U=e_m^T .
\end{equation} 

We set
$\begin{bmatrix}
Q_{m}  \\ P_{m}
\end{bmatrix} U=\begin{bmatrix}
Q'_{m}  \\  P'_{m}
\end{bmatrix} $ and $U^TB_mU=T'_{m}$, where $T'_{m}$ is a new upper Hessenberg matrix.
According to  \eqref{10}, \eqref{37or318} and \eqref{38or319},
we find the orthonormal matrix $W=W_1\times U$ and get the following equation
\begin{equation}\label{39or320}
H
\begin{bmatrix}
Q'_{m}  \\  P'_{m}
\end{bmatrix}
=
\begin{bmatrix}
Q'_{m}  \\ P'_{m}
\end{bmatrix}
T'_m+b_{m+1,m}
\begin{bmatrix}
q_{m+1} \\ p_{m+1}
\end{bmatrix}
e_{m}^{T}.
\end{equation} 
Comparing \eqref{39or320} with \eqref{10},  we complete the proof.
\end{proof}

We can get the new $k$-step GSOAR decomposition by computing the first $k$ columns of equation \eqref{39or320}
\begin{equation}\label{40or321}
H
\begin{bmatrix}
Q'_{k}  \\
P'_{k}
\end{bmatrix}
=
\begin{bmatrix}
Q'_{k}  \\
P'_{k}
\end{bmatrix}
T'_k+t'_{k+1,k}
\begin{bmatrix}
q_{k+1} \\
p_{k+1}
\end{bmatrix}
e_{k}^{T},
\end{equation} 
where the columns $\{q'_1,q'_2,\cdots ,q'_k\}$ of $Q'_{k}$ are orthonormal vectors.
\\

\begin{alg}{\bf The implicitly restarted GSOAR method with all shifts}
\begin{enumerate}
\item[1:]Run the $m$-step GSOAR procedure to generate orthonormal matrix $Q_m$,
         the columns of $Q_m$ span the subspace $\mathcal{G}_m(A,B;u_1,u_2)$;
\item[2:] Do

\hspace{0.2cm} a:  Compute the wanted approximate eigenpairs;

\hspace{0.2cm} b:  If convergence then break, else:

\hspace{0.2cm} c:  Compute $2p$ exact or refined shifts;

\hspace{0.2cm} d: Apply $2p$ implicit shifted QR algorithm to $T_m$;

\hspace{0.2cm} e: Compute $W$ as described in Theorem 3.1 and get equation \eqref{39or320}
\begin{displaymath}
H
\begin{bmatrix}
Q'_{m}  \\
P'_{m}
\end{bmatrix}
=
\begin{bmatrix}
Q'_{m}  \\
P'_{m}
\end{bmatrix}
T'_m+b_{m+1,m}
\begin{bmatrix}
q_{m+1} \\
p_{m+1}
\end{bmatrix}
e_{m}^{T};
\end{displaymath}
\hspace{0.5cm}  f: Truncate the above equation to get \eqref{40or321}
\begin{displaymath}
H
\begin{bmatrix}
Q'_{k}  \\
P'_{k}
\end{bmatrix}
=
\begin{bmatrix}
Q'_{k}  \\
P'_{k}
\end{bmatrix}
T'_k+t'_{k+1,k}
\begin{bmatrix}
q_{k+1} \\
p_{k+1}
\end{bmatrix}
e_{k}^{T};
\end{displaymath}
\hspace{0.5cm} g:  Expand  the $k$-step GSOAR decomposition to $m$.

End do
\end{enumerate}
\end{alg}

In Algorithm 3, in order to protect $k$ desired eigenvalues, we usually preserve $k+l$ instead of $k$ restarting vectors. Here, $l$ is a small nonnegative integer, such as 3 \cite{sorensen}. One restarting in Algorithm 3 needs $m-1$ times Householder transformations and $3m-6$ times matrix multiplications using $O(m^3)$ flops totally more than traditional implicitly restarted methods. Nevertheless, the dimensions of operated vectors and matrices are reducing as the procedure runs. We use all the shift candidates to improve the efficiency of each restarting, which ultimately improves the overall efficiency.

\section{Numerical experiments}

Several numerical experiments are presented in this section to demonstrate the practicability and efficiency of Algorithm 3.
We use IGSOAR to indicate the implicitly restarted second-order Arnoldi method with all exact shifts and IRGSOAR is the implicitly restarted second-order Arnoldi method with all refined shifts.
We will show the superiority of the new methods. In addition, we will make a comparison between them and the corresponding counterparts IGSOAR0, IRGSOAR0 proposed in \cite{jiasun}. All the examples are run in the same environment: Win7, 64-bit operating system, Intel(R) Core(TM) i5-2430M CPU 2.4GHz, RAM 4GB using Matlab R2010a.

For all examples, nonnegative integer $m$ stands for the dimension of projection subspace, $k$ is the number of desired eigenvalues, $f$ is the dimension of castoff subspace. So there are $p$ shifts $(p=2f)$. In all tables, $tol$ denotes the convergent criterion, TOTAL expresses the total CPU time. SOAR stands for the CPU time of computing the projected QEP and expanding GSOAR procedure from step $k+1$ to step $m$, RESTART is the CPU time of performing implicit restarting and the time of solving projected QEP is denoted by FIND, unit is second. The abscissa axis is restart times and the vertical axis is the relative residual norms in following figures.

\begin{example}
This example is tested in \cite{huangwq,jiasun}. The $8010\times 8010$ matrices are:
\begin{displaymath}
M=-4\pi ^2h^2I_{q-1}\otimes (I_q-\frac 12e_qe_q^T)\ \ ,\ \ C=2\pi i\frac h\xi I_{q-1}\otimes (e_qe_q^T),
\end{displaymath}
\begin{displaymath}
K=I_{q-1}\otimes D_q+T_{q-1}\otimes(-I_q+\frac 12e_qe_q^T),
\end{displaymath}
where $h=1/90$ is the mesh size, $q=1/h$, $\otimes$ is the Kronecker product, $\xi$ stands for impedance, $D_q=tridiag(-1,4,-1)-2e_qe_q^T$, $T_{q-1}=tridiag(1,0,1)$.  We adopt IRGSOAR, IGSOAR, IRGSOAR0 and IGSOAR0 to compute 6 desired approximate eigenvalues nearest to the origin with $m=12$, $f=5$, respectively. Table 1 and Figure 2 show the results.

\begin{table}[!h]
\begin{center}
\caption{Example 4.1, $tol={10}^{-10}$}
\label{table1}
\begin{tabular}[l]{@{}lccccccc}\hline
  Algorithm&$m$&$f$&restarts&TOTAL&SOAR&RESTART&FIND\\\hline
  IRGSOAR&12&5&3&0.600&0.357&0.086&0.147 \\
  IGSOAR&12&5&3&0.521&0.324&0.076&0.109\\
  IRGSOAR0&12&5&5&1.091&0.561&0.187&0.294\\
  IGSOAR0&12&5&5&1.017&0.571&0.176&0.211 \\\hline
\end{tabular}
\end{center}
\end{table}

\begin{figure}[!h]
\begin{center}
\resizebox*{5cm}{!}{\includegraphics{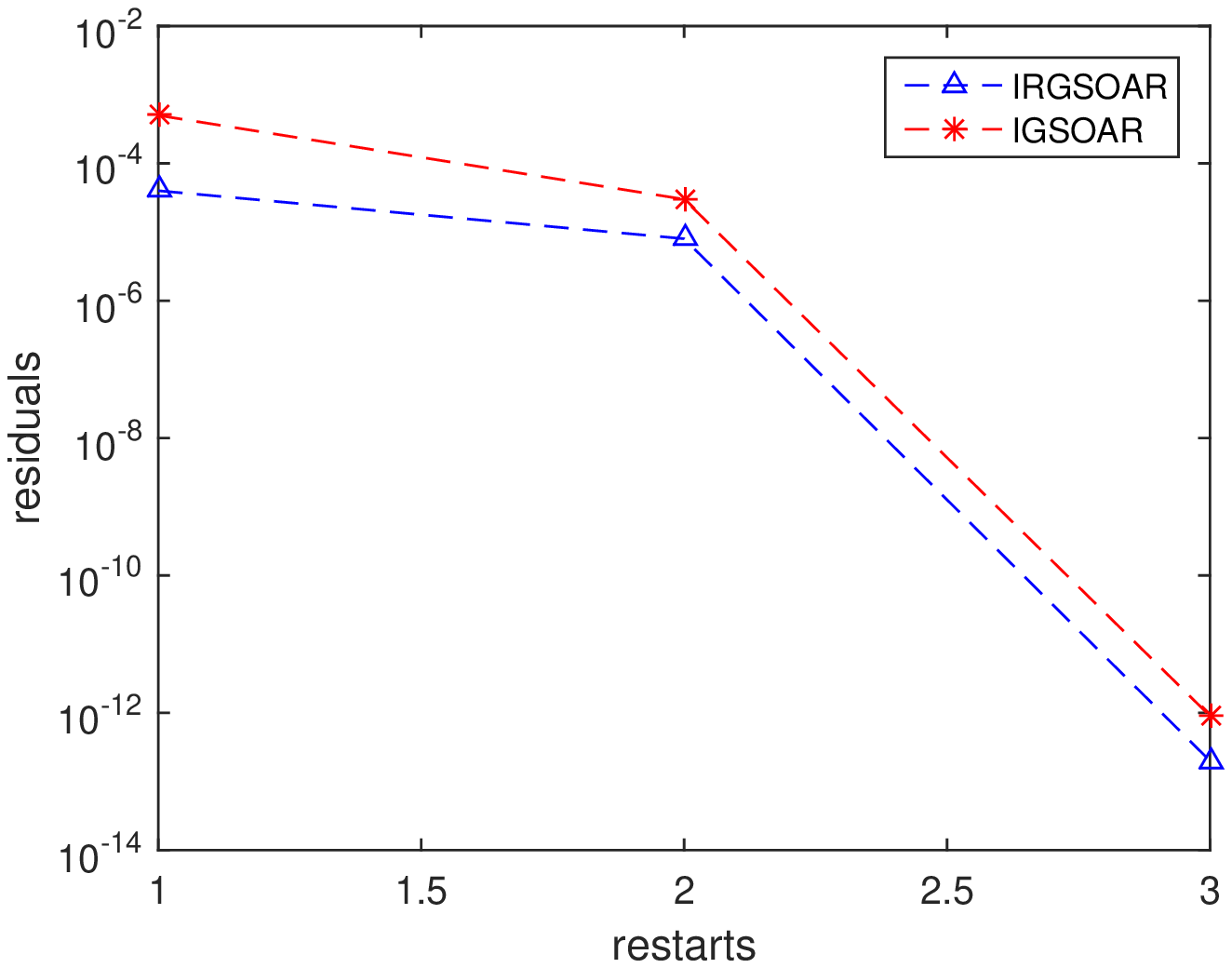}}\hspace{5pt}
\resizebox*{5cm}{!}{\includegraphics{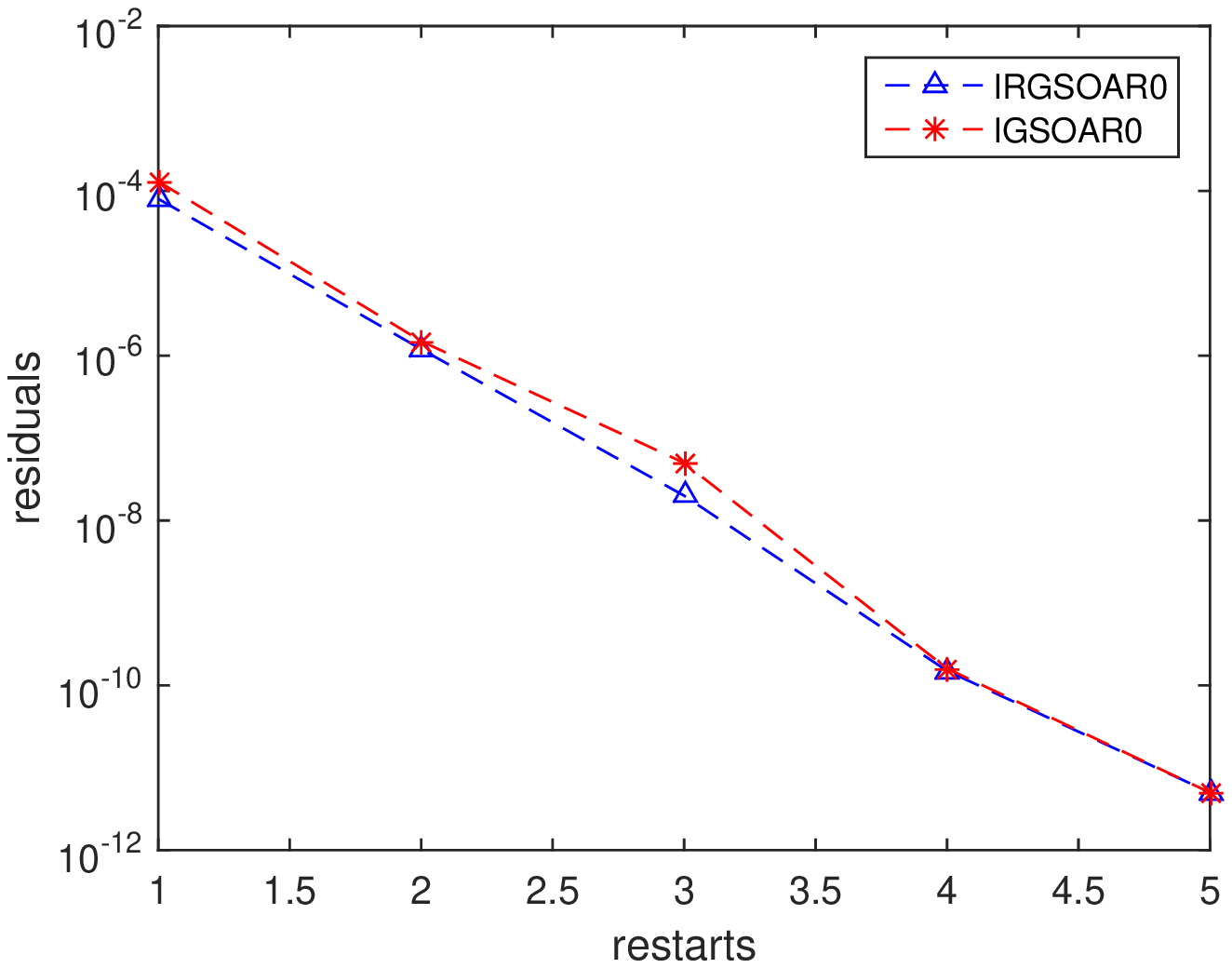}}
\caption{Example 4.1, $m=12$, $f=5$}
\label{pic2}
\end{center}
\end{figure}

It can be found from Table 1 and Figure 2 that IRGSOAR and IGSOAR work very well, their convergent speed is fast, moreover, total time is much less than IRGSOAR0 and IGSOAR0. New algorithms just need 3 restarts in contrast with 5 times of previous algorithms with part of shift candidates. The IRGSOAR and IGSOAR methods use 114 householder transformations with 0.012s. However, total time is saved much more than 0.012s. We notice that IRGSOAR is a little more efficient than the other¡¯s.
\end{example}

\begin{example}
{\bf(a)} This example is from \cite{jiasun} (cf. Example 4). The $5000\times 5000$ matrices are:
\begin{displaymath}
M=I\ ,\ C=\tau \cdot tridiag(-1,3,-1)\ ,\ K=\kappa \cdot tridiag(-1,3,-1),
\end{displaymath}
\noindent where $\tau$ =10, $\kappa$ =5. We adopt IRGSOAR, IGSOAR, IRGSOAR0 and IGSOAR0 to compute 6 desired approximate eigenvalues nearest to target $\sigma =-13+0.4i$ with $m=40$, $f=28$, respectively. Table 2 and Figure 3 show the results.

\begin{table}[!h]
\begin{center}
\caption{Example 4.2.(a), $tol={10}^{-10}$}
\label{table2}
\begin{tabular}[l]{@{}lccccccc}\hline
  Algorithm&$m$&$f$&restarts&TOTAL&SOAR&RESTART&FIND\\\hline
  IRGSOAR&40&28&4&1.803&0.558&0.644&0.581 \\
 IGSOAR&40&28&6&2.567&0.700&1.174&0.673\\
  IRGSOAR0&40&28&54&27.390&7.955&10.008&8.590 \\
  IGSOAR0&40&28&65&29.024&9.617&12.239&6.173\\\hline
\end{tabular}
\end{center}
\end{table}

\begin{figure}
\begin{center}
\resizebox*{5cm}{!}{\includegraphics{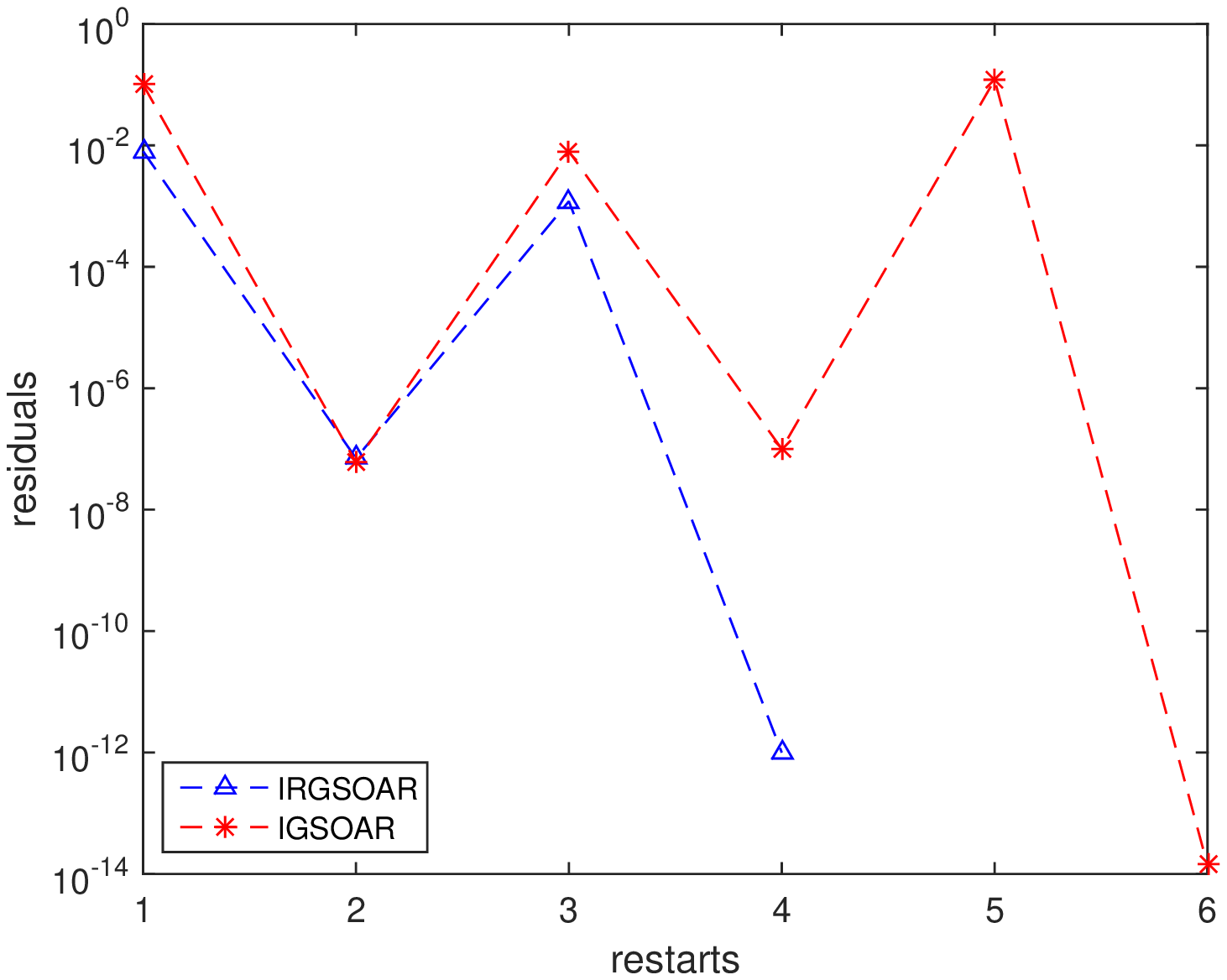}}\hspace{5pt}
\resizebox*{5cm}{!}{\includegraphics{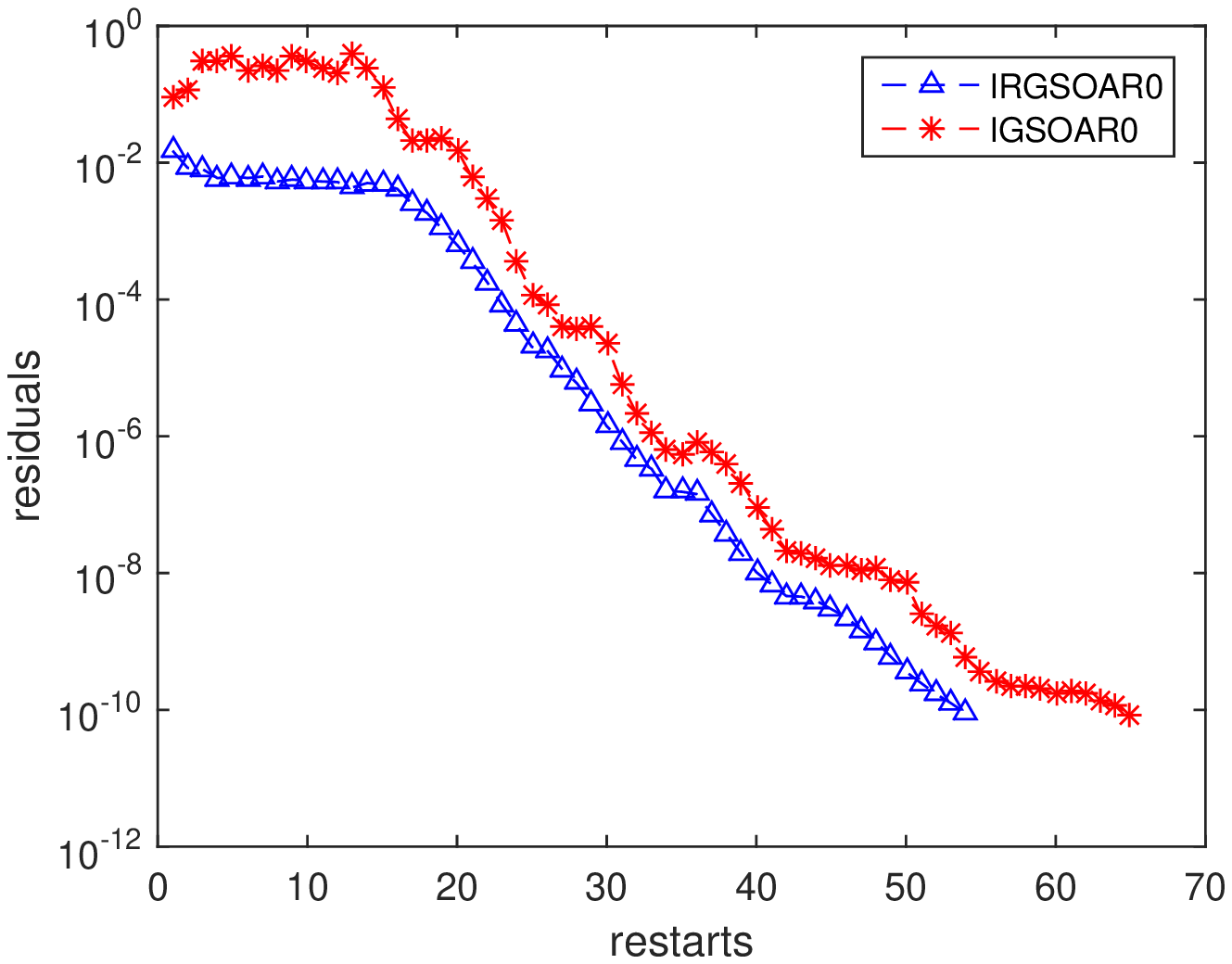}}
\caption{Example 4.2.(a), $m=40$, $f=28$}
\label{pic3}
\end{center}
\end{figure}

{\bf(b)} If $f=30$, other conditions are the same. Run the methods respectively. Table 3 and Figure 4 show the results.

\begin{table}[!h]
\begin{center}
\caption{Example 4.2.(b), $tol={10}^{-10}$}
\label{table3}
\begin{tabular}[l]{@{}lccccccc}\hline
  Algorithm&$m$&$f$&restarts&TOTAL&SOAR&RESTART&FIND\\\hline
  IRGSOAR&40&30&5&2.341&0.616&0.979&0.730 \\
  IGSOAR&40&30&6&2.580&0.81&1.177&0.572\\
  IRGSOAR0&40&30&54&26.880&8.154&10.013&7.877 \\
  IGSOAR0&40&30&59&25.508&8.514&10.665&5.416\\\hline
\end{tabular}
\end{center}
\end{table}

\begin{figure}
\begin{center}
\resizebox*{5cm}{!}{\includegraphics{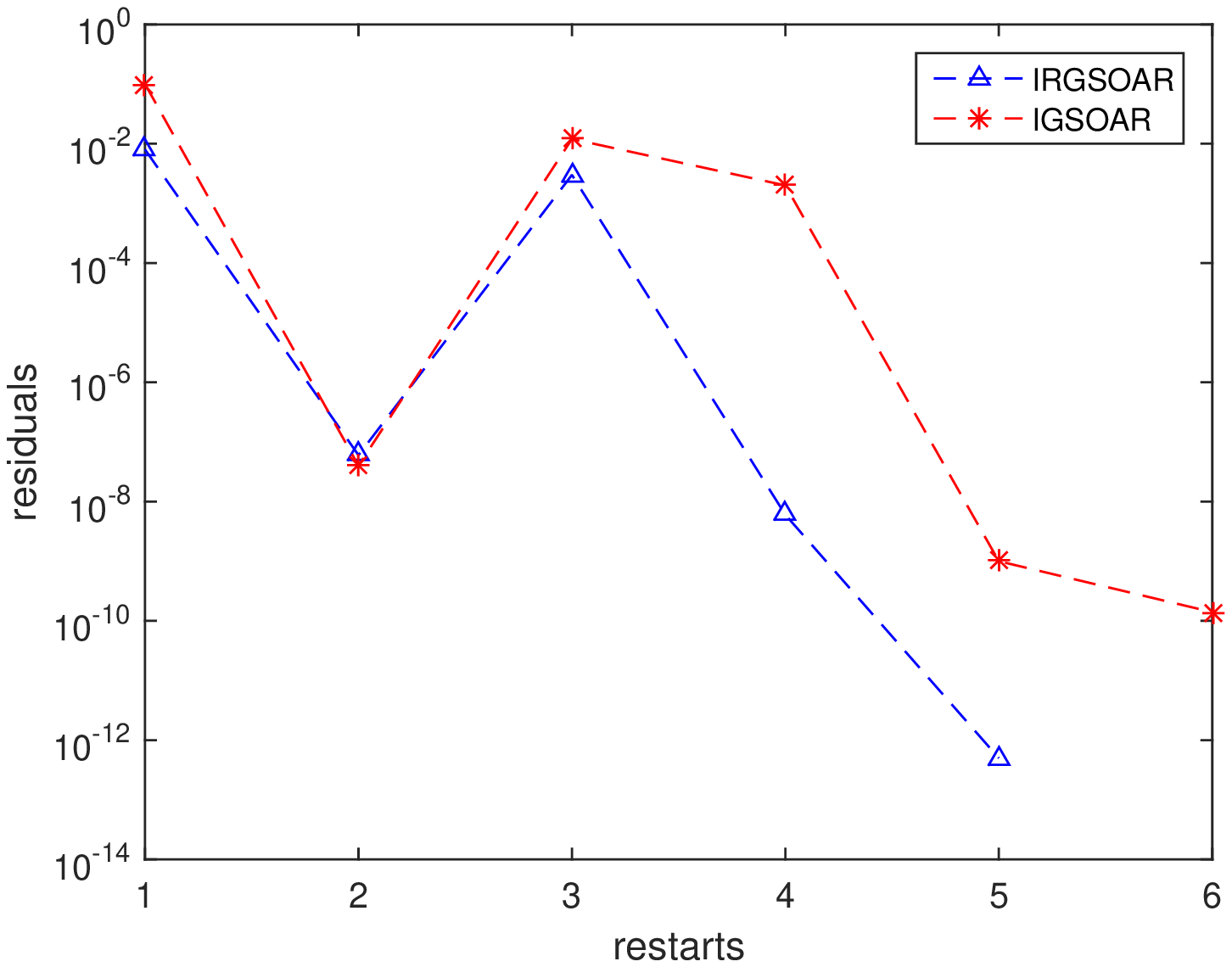}}\hspace{5pt}
\resizebox*{5cm}{!}{\includegraphics{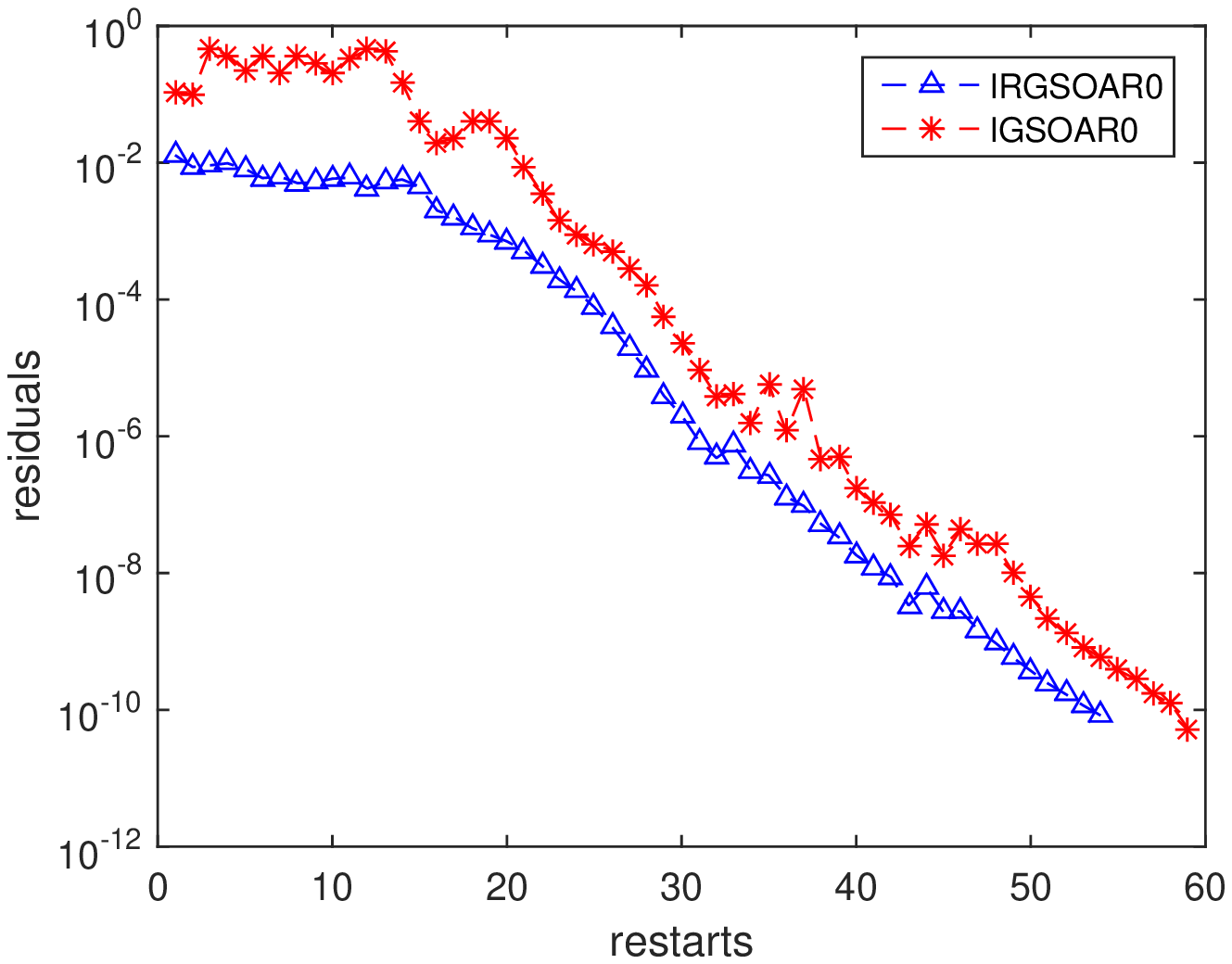}}
\caption{Example 4.2.(b), $m=40$, $f=30$}
\label{pic4}
\end{center}
\end{figure}

We can see the marked superiority of IRGSOAR and IGSOAR algorithms. When $f=28$, previous algorithms restarted more than 50 times to reach the convergent requirement. While IRGSOAR and IGSOAR algorithms just need 4 and 6 restarts, respectively. In addition, IRGSOAR saved almost 26s, in contrast to consuming 0.028s for Householder transformation. When $f=30$, previous algorithms still needed more than 50 restarts to reach the convergent requirement. While IRGSOAR and IGSOAR algorithms just need 5 and 6 restarts, respectively. IRGSOAR always need less restarts and less running time even when $f$ changes. Example 2 tells us the implicit restarting methods are more advantageous in total running time and restart times from shifts increasing.
\end{example}

\begin{example}
{\bf(a)} The size of this problem is $n=5000$ and the matrices are \\
\begin{displaymath}
M=I_n\ , \
C=
\begin{bmatrix}
8&-4&0&  \\
2&12&\ddots &0\\
0&\ddots &12&-4\\
 &0&2&8
\end{bmatrix}\ , \
K=
\begin{bmatrix}
2&2&0&  \\
-1&3&\ddots &0\\
0&\ddots &3&2\\
 &0&-1&2
\end{bmatrix}.
\end{displaymath}
We adopt IRGSOAR, IGSOAR, IRGSOAR0 and IGSOAR0 to compute 6 desired approximate eigenvalues nearest to target $\sigma =-10-0.8i$ with $m=26$, $f=15$, respectively. Table 4 and Figure 5 show the results.

\begin{table}[!h]
\begin{center}
\caption{Example 4.3.(a), $tol={10}^{-10}$}
\label{table4}
\begin{tabular}[l]{@{}lccccccc}\hline
  Algorithm&$m$&$f$&restarts&TOTAL&SOAR&RESTART&FIND\\\hline
  IRGSOAR&26&15&2&1.014&0.102&0.101&0.391 \\
  IGSOAR&26&15&3&0.936&0.194&0.175&0.129\\
  IRGSOAR0&26&15&100&53.711&6.523&12.169&8.312\\
  IGSOAR0&26&15&100&47.8455&6.494&12.095&5.497 \\\hline
\end{tabular}
\end{center}
\end{table}

\begin{figure}[!h]
\begin{center}
\resizebox*{5cm}{!}{\includegraphics{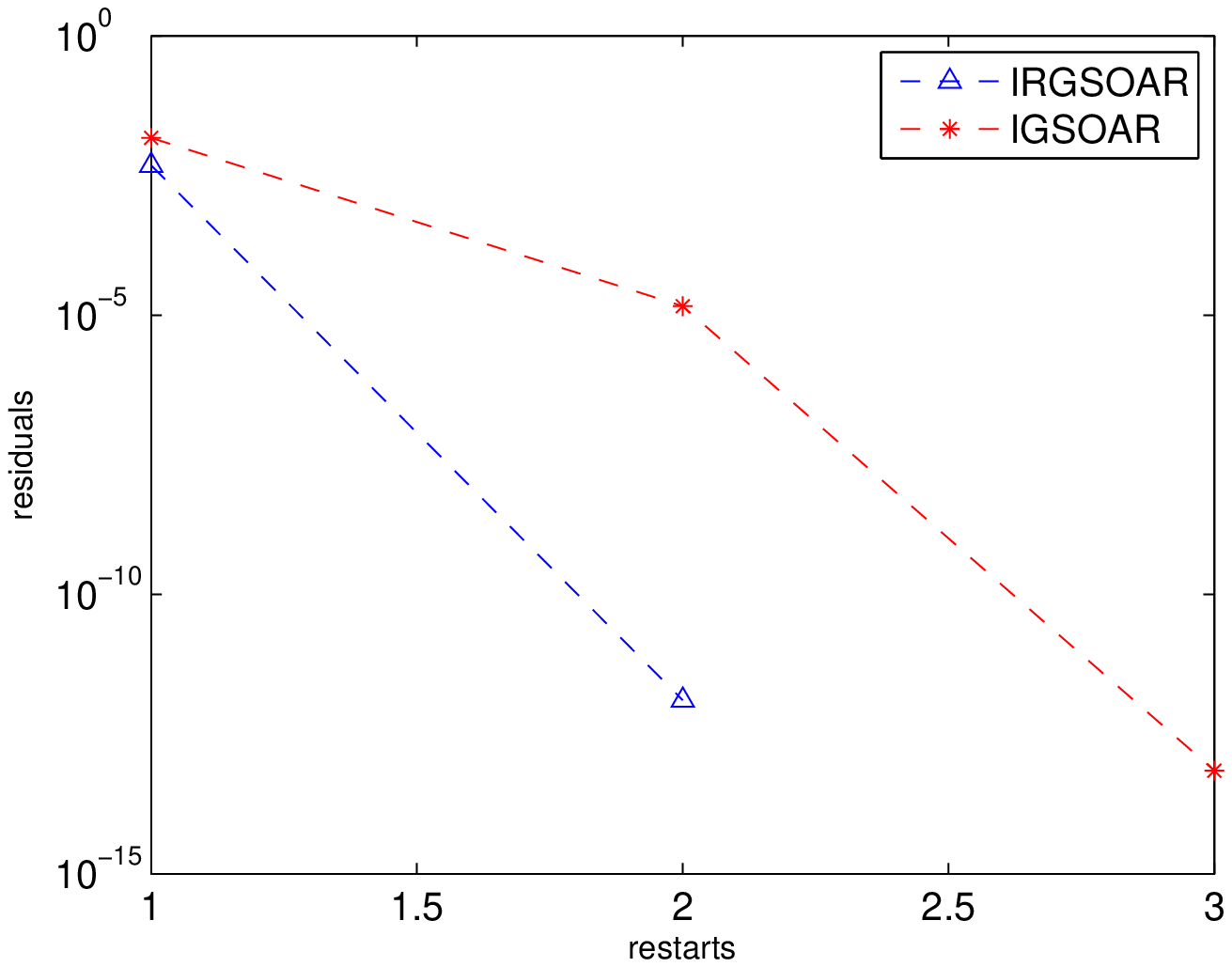}}\hspace{5pt}
\resizebox*{5cm}{!}{\includegraphics{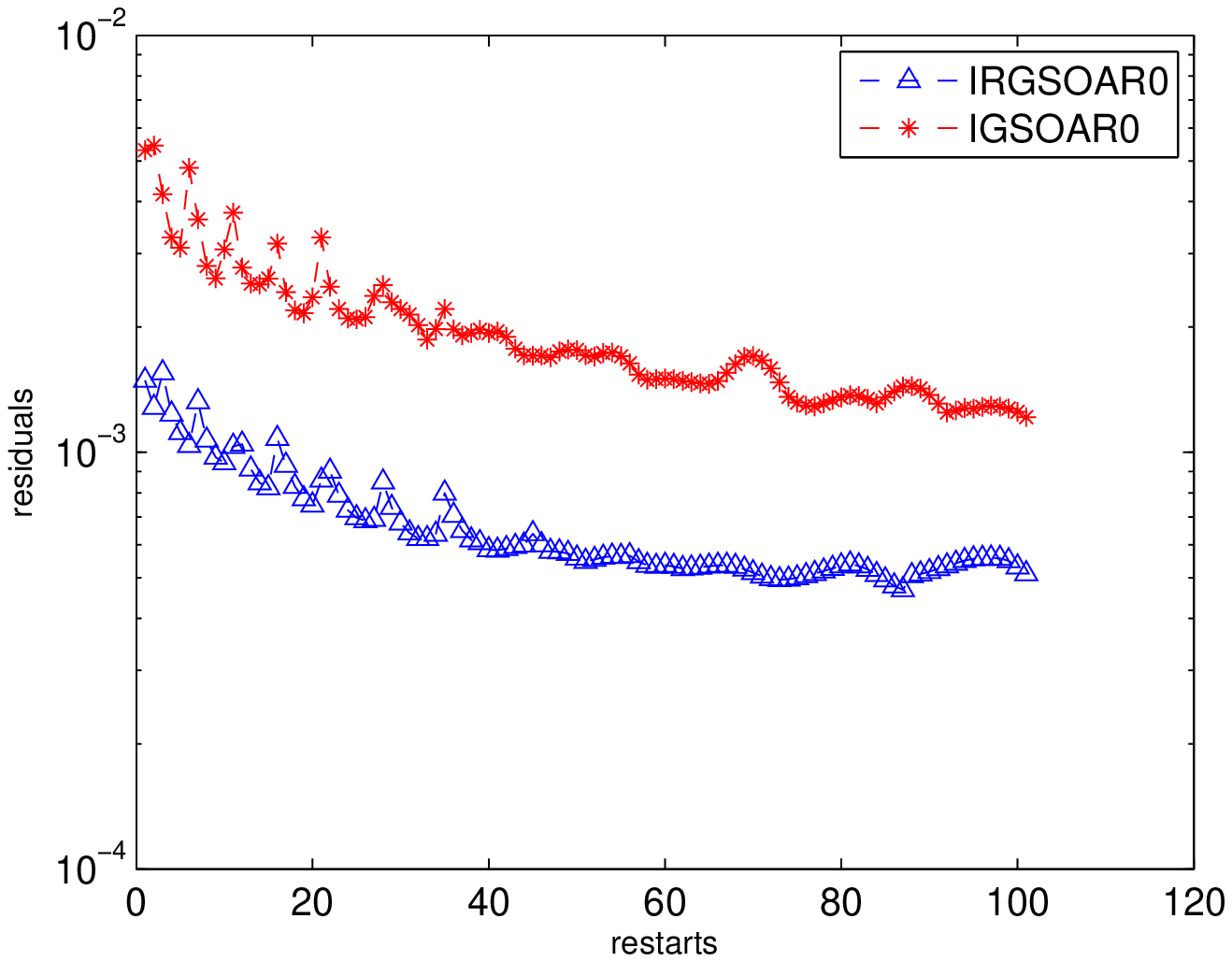}}
\caption{Example 4.3.(a), $m=26$, $f=15$}
\label{pic5}
\end{center}
\end{figure}

{\bf(b)} If $f=13$, other conditions are the same. Run the methods respectively. Table 5 and Figure 6 show the results.

\begin{table}[!h]
\begin{center}
\caption{Example 4.3.(b), $tol={10}^{-10}$}
\label{table5}
\begin{tabular}[l]{@{}lccccccc}\hline
  Algorithm&$m$&$f$&restarts&TOTAL&SOAR&RESTART&FIND\\\hline
  IRGSOAR&26&13&2&0.858&0.186&0.109&0.166 \\
  IGSOAR&26&13&3&0.858&0.186&0.171&0.158\\
  IRGSOAR0&26&13&100&50.513&5.500&11.409&7.981\\
  IGSOAR0&26&13&100&46.067&5.710&11.862&5.283 \\\hline
\end{tabular}
\end{center}
\end{table}

\begin{figure}[!h]
\begin{center}
\resizebox*{5cm}{!}{\includegraphics{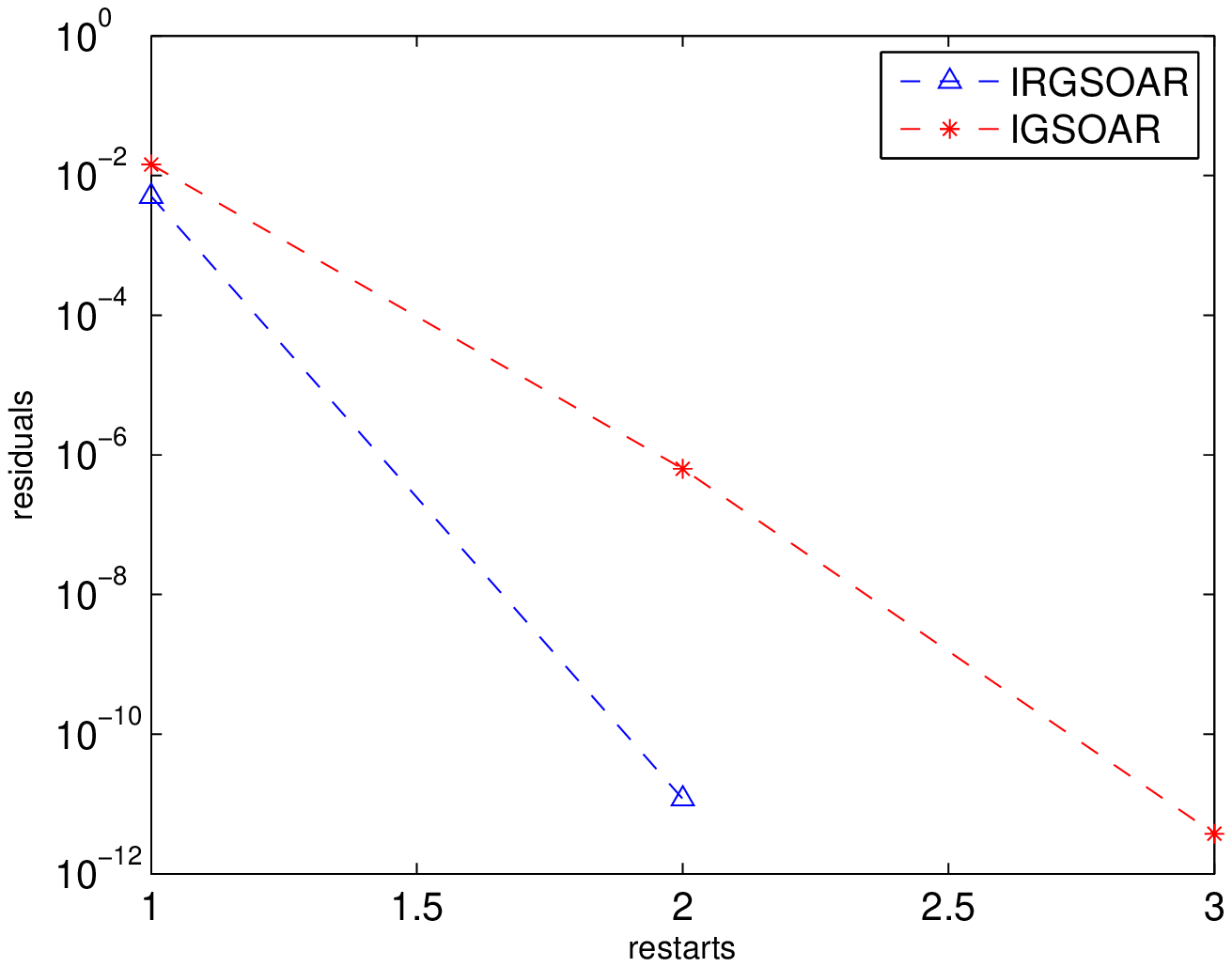}}\hspace{5pt}
\resizebox*{5cm}{!}{\includegraphics{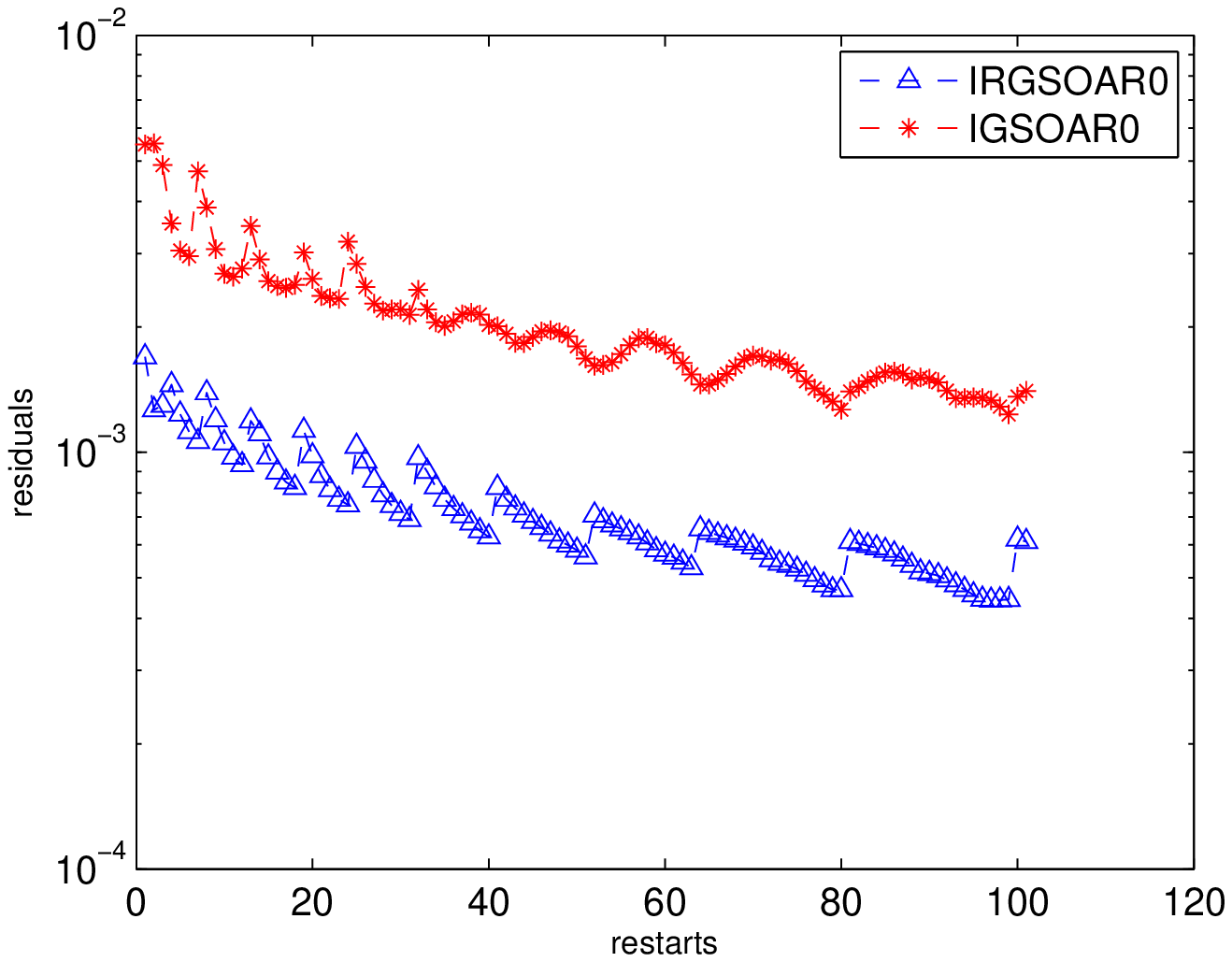}}
\caption{Example 4.3.(b), $m=26$, $f=13$}
\label{pic6}
\end{center}
\end{figure}

It can be seen from the results that the previous algorithms do not satisfy the convergent criterion yet after 100 restarts. However, the algorithms proposed in our work perform well on this problem. They only need two or three restarts to solve the large QEP.
\end{example}

We have developed implicitly restarted algorithms
with certain exact and refined shifts proposed for generalized second-order Arnoldi method.
Unlike the linear eigenvalue problem, for the QEP,
there are more shift candidates than the shifts allowed. To overcome this problem,
we have presented an efficient and reliable algorithm for using all shift
candidates. Numerical experiments have demonstrated that using more shifts
to improve each restart can boost efficiency significantly.

\Acknowledgements{This work was supported by National Natural Science Foundation of China (Grant No. 11201020).}



\begin{thebibliography}{99}
\bahao\baselineskip 11.5pt

\bibitem{baidemmel}
Bai Z, Demmel J, Dongarra J,Ruhe A and van der Vorst H A. Templates for the Solution of Algebraic Eigenvalue Problems: A Practical Guide. SIAM, Philadelphia, PA, 2000.

\bibitem{baisu}
Bai Z and Su Y. SOAR: A second-order Arnoldi
method for the solution of the quadratic eigenvalue problem.
SIAM J. Matrix Anal, 2003, 26: 640--659.

\bibitem{betcke}
Betcke T, Higham N J, Mehrmann V, Schr\"{o}der C and Tisseur F. NLEVP: A collection of nonlinear eigenvalue problems. users' guide, MIMS EPrint, 2010.

\bibitem{huangt}
Huang H M, Jia Z and Lin W W. On the convergence of Ritz pairs and refined Ritz vectors for quadratic eigenvalue problems. BIT Numerical Mathematics, 2013, 53: 941--958.

\bibitem{huangwq}
Huang W Q, Li T, Li Y T and Lin W W. A semiorthogonal generalized Arnoldi method and its variations for quadratic eigenvalue problems. Numerical Linear Algebra with Applications, 2013, 20: 259--280.

\bibitem{jiarefined}
Jia Z. Refined iterative algorithms based on Arnoldi process for large unsymmetric eigenproblems. Linear Algebra Appl, 1997, 259: 1--23.

\bibitem{jia}
Jia Z. Polynomial characterizations of the approximate eigenvectors by the refined Arnoldi method and an implicitly restarted refined Arnoldi algorithm. Linear algebra and its applications, 1999, 287: 191--214.

\bibitem{jia02}
Jia Z. The refined harmonic Arnoldi method
and an implicitly restarted refined algorithm for computing interior
eigenpairs of large matrices. Appl. Numer. Math, 2002, 42: 489--512.

\bibitem{jiasun}
Jia Z and Sun Y. Implicitly restarted generalized second-order Arnoldi type algorithms for the quadratic eigenvalue problem. Taiwanese Journal of Mathematics, 2015, 19: 1--30.

\bibitem{meerbergen}
Meerbergen K. The quadratic Arnoldi method for the solution of the quadratic eigenvalue problem. SIAM Journal on Matrix Analysis and Applications, 2008, 30: 1463--1482.

\bibitem{otto}
Otto C. Arnoldi and Jacobi¨CDavidson methods for quadratic eigenvalue problems. diploma thesis, Institut f\"{u}r Mathematik,
Technische Universit\"{a}t Berlin, Germany, 2004.

\bibitem{ruhe}
Ruhe A. Algorithms for the nonlinear eigenvalue problem. SIAM Journal on Numerical Analysis, 1973, 10: 674--689.

\bibitem{sorensen}
Sorensen D C. Implicit application of polynomial filters in a k-step Arnoldi method. SIAM Journal on Matrix Analysis and Applications, 1992, 13: 357--385.

\bibitem{stewart}
Stewart G W. A Krylov--Schur Algorithm for Large Eigenproblems. SIAM Journal on Matrix Analysis and Applications, 2002, 23: 601--614.

\bibitem{tisseur}
Tisseur F and Meerbergen K. The quadratic eigenvalue problem. SIAM Rev, 2001, 43: 235-286.

\bibitem{zhou}
Zhou L, Bao L, Lin Y, Wei Y and Wu Q. Restarted generalized Krylov subspace methods for solving quadratic eigenvalue problems. Inter. J. Comput. Math. Sci, 2000, 4: 148--155.


\end{thebibliography}
\end{document}